\documentclass[11pt]{amsart}

\usepackage[colorlinks=true,urlcolor=black,anchorcolor=black,citecolor=black,linkcolor=black]{hyperref}
\usepackage{datetime}
\usepackage{enumerate
}\usepackage{amsmath}

\DeclareMathOperator{\KFD}{\Delta}

\usepackage{amssymb, amsfonts, amsmath, amsthm}
\usepackage{tabularx}
\usepackage{xcolor} 
\usepackage{graphicx}
\usepackage{array}      
\usepackage{booktabs}
\usepackage{mathrsfs}
\usepackage{subcaption}
 
\usepackage{graphicx}
\newtheorem{theorem}{Theorem}[section]
\newtheorem{lemma}[theorem]{Lemma}
\newtheorem{proposition}[theorem]{Proposition}
\newtheorem{corollary}[theorem]{Corollary}
\newtheorem{prop-and-def}[theorem]{Proposition and Definition}
\newtheorem{proposition-and-definition}[theorem]{Proposition and Definition}

\theoremstyle{definition}

\newtheorem{notation}[theorem]{Notation}
\newtheorem{remark}[theorem]{Remark}
\newtheorem{definition-and-remark}[theorem]{Definition and Remark}
\newtheorem{remark-and-notation}[theorem]{Remark and Notation}
\newtheorem{remark-and-definition}[theorem]{Remark and Definition}
\newtheorem{notation-and-remark}[theorem]{Notation and Remark}

\newtheorem{assumption}[theorem]{Assumption}

\numberwithin{equation}{section}

\renewcommand{\Im}{\mathop{\mathrm{Im}}}
\renewcommand{\Re}{\mathop{\mathrm{Re}}}
\newcommand{\bd}{\mathop{\mathrm{bd}}}
\newcommand{\spec}{\mathop{\mathrm{spec}}}

\newcommand\bub[1]{\mathring{#1}}
\newcommand{\norm}[1]{\lVert#1\rVert}
\newcommand{\abs}[1]{\left\lvert #1 \right\rvert}  

\newcommand{\cA}{ \mathcal{A} }
\newcommand{\cB}{ \mathcal{B} }
\newcommand{\cD}{ \mathscr{D} }

\newcommand{\cM}{ \mathscr{M} }

\newcommand{\bC}{ \mathbb{C} }

\newcommand{\bH}{ \mathbb{H} }
\newcommand{\bN}{ \mathbb{N} }

\newcommand{\bR}{ \mathbb{R} }

\newcommand{\ee}{ \varepsilon }
\newcommand{\qi}{ \mathrm{i} }

\newcommand{\pp}{p}

\title[On the Brown measure of $X + \qi Y$, with $X,Y$ selfadjoint 
and $Y$ free Poisson]{On the Brown measure of $x+ \qi y$, 
with $x,y$ selfadjoint \\
and $y$ free Poisson} 

\author[F. Lehner]{Franz Lehner}
\address{Franz Lehner: Institut f\"ur Diskrete Mathematik,
Technische Universit\"at Graz,
Steyrergasse 30,
A-8010 Graz, Austria}
\thanks{FL: This research was funded in part by the Austrian Science Fund (FWF) grant [10.55776/I6232] BOOMER (WEAVE).
}

\email{lehner@math.tugraz.at}

\author[A. Nica]{Alexandru Nica}
\address{Alexandru Nica: Department of Pure Mathematics, 
	University of Waterloo, Ontario, Canada.}
\email{anica@uwaterloo.ca}
\thanks{AN: research supported by a Discovery Grant from NSERC, Canada.}

\author[K. Szpojankowski]{Kamil Szpojankowski}
\address{ Kamil Szpojankowski: Institute of Mathematics of the Polish Academy of Sciences, ul. \'Sniadeckich 8, 00-656 Warszawa, Poland, and 
	Faculty of Mathematics and Information Science,
	Warsaw University of Technology, Koszykowa 75, 00-662 Warszawa, Poland.} 
\email{kszpojankowski@impan.pl, kamil.szpojankowski@pw.edu.pl  }
\thanks{KSz: This research was funded in part by National Science Centre, Poland WEAVE-UNISONO grant BOOMER 2022/04/Y/ST1/00008.}

\author[P. Zhong]{Ping Zhong}
\address{ Ping Zhong: Department of Mathematics, University of Houston, Houston, TX 77204-3008, USA } 
\email{pzhong@central.uh.edu }
\thanks{PZ: supported in part by NSF grant LEAPS-MPS-2516951 and NSF CAREER Award DMS-2516950}

\date{ \today, \currenttime UTC}

\begin{document}

\begin{abstract}  
Let $x,y$ be freely independent selfadjoint elements in a $W^{*}$-probability space, 
where $y$ has free Poisson distribution of parameter $\pp$.  

We pursue a methodology 
for computing the Brown measure of $x + \qi y$, 
which relies on the matrix-valued subordination function $\Omega$ of the hermitization 
of $x + \qi y$, and on the fact that $\Omega$ has an explicitly described left inverse $H$. 
Our main point is that the Brown measure of $x + \qi y$ becomes more approachable when 
it is reparametrized via a certain change of variable $h : \cD \to \cM$, with $\cD, \cM$
open subsets of $\bC$, where $\cD$ and $h$ are defined in terms of the aforementioned 
left inverse $H$, and $\mathrm{cl} \,(\cM)$ contains the support of the absolutely continuous part of Brown measure. 
More precisely, we find (with some conditions on the distribution of $x$) the following formula:
\[
f(s + \qi \, t) =\frac{1}{2 \pi}\left[\frac{1}{t}\left(\frac{\partial \alpha}{\partial s}
+\frac{\partial \beta}{\partial t}\right)-\frac{1}{t}-\frac{\beta}{t^2}\right],
\ \ s + \qi \, t \in \cM,
\]
where $f$ is the density of the absolutely continuous part of the Brown measure 
and the functions $\alpha, \beta : \cM \to \bR$ are the real and respectively the 
imaginary part of $h^{-1}$.

We show that if $x$ has an atom  $\alpha$ with $\mu_x(\alpha)>p$,
then the Brown measure of $x+\qi y$ has an atom of mass $\mu_x(\alpha)-p$ at the same point $\alpha$. Moreover we prove that if $\alpha_1,\ldots,\alpha_k$  is the list of atoms of $x$ with mass bigger than $p$,
then the Brown measure of $x+\qi y$ is supported on $\mathrm{cl}(\cM)\cup\{\alpha_1,\ldots,\alpha_k\}$.
\end{abstract}

\maketitle

\tableofcontents

\section{Introduction}

The notion of Brown measure, introduced by L.G.~Brown in \cite{Brown}, serves as 
a replacement for the notion of spectral distribution measure for non-normal elements
in a $W^*$-probability space.  The study of this notion became 
popular starting with the early 2000's, when it was found that Brown measures can 
be explicitly computed in special cases arising in free probabilisty. 
In some of these special cases, the Brown measure was also found to be the correct
predictor for the limit eigenvalue distribution of the corresponding non-Hermitian 
matrix models.  A more detailed review of how the Brown measure is defined and some 
references to the relevant research literature are given in Section 2.1 below.

Explicit computations of Brown measures are known to be challenging, but the research
literature of the last few years has created a gallery of interesting examples where 
such explicit computations can be carried through.  
There are two main 
lines of approach that have been developed towards this purpose.  One of these lines 
(cf.~\cite{BelinschiSniadySpeicher,BelinschiYinZhong2024, BercoviciZhong2022, HeltonSpeicherMai, Zhong2021Brown}) 
relies on matrix-valued constructions encoding the non-normal element whose Brown measure 
is being studied, followed by an analysis of subordination functions for some matrix-valued 
Cauchy transforms which appear in this way.  The other line 
(cf.~\cite{DHK2022Brown, HallHo2022LMP, HallHo2023ptrf, HoZhong2023Brown, DemniHamdi2022jfa}) 
pursues an analysis of partial differential equations which arise in connection with free additive 
and free multiplicative Brownian motions.

To put things into perspective, we mention that the first explicit computations of
Brown measures in free probabilistic setting were in connection to products $x \cdot y$
where $x,y$ are selfadjoint and freely independent cf.~\cite{HaagerupLarsen,BianeLehner2001}.  
For same kind of $x$ and $y$, Brown measures of elements $x + \qi y$ were also tackled in
subsequent work -- we mention here the work in \cite{BelinschiYinZhong2024, BCC2014cpam, Zhong2021Brown}, 
where $y$ is the imaginary part of a circular (or more generally elliptic) element free 
from $x$, and the one in \cite{HallHo2022LMP}, where $y$ is semicircular. 

The present paper, too, considers elements of the form $x + \qi y$ where $x$ and $y$ 
are selfadjoint and freely independent.  Throughout the paper we assume that 
$y$ has free Poisson distribution with parameter $\pp$, i.e., all free cumulants of 
$y$ are equal to $\pp$.  For such $x$ and $y$,  
our main point is that the Brown measure of $x + \qi y$ is better accessible via a certain reparametrization $h : \cD \to \cM$, with $\cD, \cM$ open subsets 
of $\bC$ such that the closure $\mathrm{cl} \,(\cM)$ contains the support of  the absolutely continuous part of the Brown measure, and 
where $\cD, h$ are defined by using subordination. Moreover the atoms of the Brown measure can also be described in terms of $\cD$ and $h$. This approach is primarily following 
the first of 
the two lines mentioned above, via subordination, combined with some ideas inspired by the PDE approach in \cite{HallHo2022LMP}.

The next subsections present what one might call ``a general methodology'' underlying 
our approach.  While our discussion is focused on the case when $y$ has free Poisson 
distribution, we hope that this general methodology will apply
to more examples of elements of the form $x+\qi y$, including (on a didactical level) the case when $y$ has semicircular distribution.

In Subsections~\ref{subsection:11}-\ref{subsection:13} below we outline the main steps 
of our general methodology.  In these subsections and throughout the paper we use the following framework:

\begin{itemize} 
\item A $W^{*}$-probability space $( \cA , \varphi )$ where $\varphi$ is a faithful trace.

\item Two selfadjoint elements $x,y \in \cA$ which are freely independent and where $y$ 
has free Poisson distribution with parameter $\pp \in (0, \infty )$.  We let $\mu_x$ denote
the distribution of $x$ with respect to $\varphi$ (a compactly supported measure on $\bR$).
It is well-known that, due to the faithfulness of $\varphi$, the support of $\mu_x$ is 
equal to $\spec (x)$, the spectrum of $x$ in $\cA$.

\item We let $\mu_{x + \qi y}$ denote the Brown measure of $x + \qi y$ 
with respect to $\varphi$ (a compactly supported measure on $\bC$).
\end{itemize} 

\noindent
The main result of the paper, concerning the Brown measure $\mu_{x + \qi y}$, is
Theorem~\ref{thm:18} stated in Section~\ref{subsection:13}.  

\medskip

\subsection{Hermitization, subordination, and the function $H$.}
\label{subsection:11}
\phantom{x}

\noindent
We consider the $W^{*}$-algebra $M_2 ( \cA )$ and the selfadjoint elements 
$X,Y \in M_2 ( \cA )$ defined as
\begin{equation}   \label{eqn:1b}
X = \begin{bmatrix}  0 & x  \\ x & 0  \end{bmatrix}, 
\mbox{ and } Y = \begin{bmatrix} 0 & \qi y \\ - \qi y & 0\end{bmatrix}. 
\end{equation}
The matrix 
\[
X + Y =  
\begin{bmatrix}
    0   & x + \qi y \\  ( x + \qi y)^{*} & 0
\end{bmatrix} 
\]
is known as the {\em hermitization} of the element $x + \qi y \in \cA$.

An easily verified basic fact in matrix-valued free 
probability is that when one considers the conditional expectation 
$E : M_2(\cA)\to M_2(\bC)$ defined by
\begin{align*}
        E (A) = 
        \begin{bmatrix}
            \varphi(a_{11})&\varphi(a_{12})\\ \varphi(a_{21})&\varphi(a_{22})
        \end{bmatrix}, \mbox{ for } 
        A = \begin{bmatrix}
            a_{11}&a_{12}\\a_{21}&a_{22}
        \end{bmatrix} \in M_2 ( \cA ),
\end{align*}
the free independence of $x$ and $y$ in $( \cA, \varphi )$ entails that 
$X,Y$ from (\ref{eqn:1b}) are freely independent in the $M_2 ( \bC )$-valued 
non-commutative probability space $( M_2 ( \cA ), E )$.  A further, non-trivial,
consequence of the latter free independence is that the matrix-valued Cauchy 
transform of $X+Y$ is subordinated to the one of $X$.  We use standard notation 
for such Cauchy transforms:
\begin{itemize}
\item[--] We denote by $\bH_2^+$ the matrix upper half plane 
\[
\bH_{2}^{+} := \{ B \in M_2 ( \bC ) \mid \Im (B) 
\text{ is positive definite and invertible}  \},
\]
\noindent
and by $\bH_{2}^{-}$ the matrix lower half plane
$\bH_{2}^{-} := \{ -B \mid B \in \bH_{2}^{+} \}$,
where as usual we denote by $\Im(B)=\frac{1}{2i}(B-B^*)$
the imaginary part of a matrix.
\vspace{6pt}

\item[--] 
For $A = A^{*} \in M_2 ( \cA )$, the Cauchy transform of $A$ is the function 
$G_A : \bH_{2}^{+} \to \bH_{2}^{-}$ defined by
$G_A (B) = E ( \, (B-A)^{-1} \, ), \ \ B \in \bH_{2}^{+}.$

\vspace{6pt}

\item[--]
Throughout the paper it will be of particular interest to evaluate our
matrix-valued Cauchy transforms at matrices $B \in \bH_2^{+}$ of the special form 
\begin{equation}   \label{B-Matrix}
B = B( \lambda , \zeta ) := 
\begin{bmatrix} \zeta & \lambda \\ \overline{\lambda} & \zeta
\end{bmatrix}, 
\ \ \mbox{ with $\lambda \in \bC$ and $\zeta \in \bC^+$.} 
\end{equation}
In relation to that, we record here that for a matrix
\begin{equation}   \label{A-Matrix}
 A = \begin{bmatrix}
      0 &a \\ a^{*} & 0  \end{bmatrix} \in M_2 ( \cA ), \mbox{ with } a \in \cA,
\end{equation}
a quick calculation with Schur complements yields the explicit formula

\begin{align}    \label{eq:entriesMatrixCauchy1}
    G_A \bigl( \, B( \lambda, \qi \delta) \, \bigr) 
= \left[  \begin{array}{cc}
          -\qi \delta \varphi\left((| \lambda -a|^2 + \delta^2)^{-1}\right)
                               & \varphi\left( ( \lambda -a)(| \lambda -a|^2 + \delta^2)^{-1} \right)  \\
                               &                                                                       \\
         \varphi\left( ( \lambda - a)^{*} (| \lambda -a |^2 + \delta^2)^{-1} \right)
                               & - \qi \delta \varphi \left((| \lambda -a |^2 + \delta^2)^{-1}\right)
\end{array}   \right] ,
\end{align}

\noindent
valid for $\lambda \in \bC$ and $\delta \in (0, \infty )$.  

In particular, setting
$a =x$ (resp.~$a = x+ \qi y$) we obtain explicit formulas for
$G_X \bigl( \, B( \lambda, \qi \delta) \, \bigr)$ and 
$G_{X+Y} \bigl( \, B( \lambda, \qi \delta) \, \bigr)$, with $X, Y$ as in (\ref{eqn:1b}).

\vspace{6pt}

\item[--]
It is known (cf.~\cite[Section 3.8]{Voiculescu2000} and 
\cite[Theorem 2.2]{BelinschiMaiSpeicher}) that the matrix-valued subordination of 
$X+Y$ holds with respect to $X$ and is 
provided by a self-map $\Omega:\bH_{2}^{+}\to \bH_{2}^{+}$ such that
\begin{equation}
\label{eq:GX+YB=GXOB}
G_{X+Y} (B) = G_X \bigl( \, \Omega(B) \, \bigr), \ \ \forall \, B \in \bH_{2}^{+}.
\end{equation}
\end{itemize}

With inspiration from ideas that previously appeared, in the selfadjoint case, 
in the work of Biane \cite{Biane}, we will take advantage of the fact that the subordination 
function $\Omega$ has a tractable left inverse.  For our purposes it is in fact convenient 
to introduce this left inverse in a self-contained way, as explained 
in the next proposition.

\medskip

\begin{proposition-and-definition}  \label{prop:11}
Given the setting introduced above,
 the subordination function $\Omega$ from \eqref{eq:GX+YB=GXOB}
has a left inverse $H \bigl( \, \Omega (B) \, \bigr) = B$, 
$\forall \, B \in \bH_{2}^{+}$
which is given by the formula
\begin{equation} \label{eqn:H_Definition}
    H(B) = B + \pp \,  \bigl( J - G_{X}(B) \bigr)^{-1},
    \ \ B \in \bH_{2}^{+},
\end{equation}
where
 $J$ is the Pauli matrix
\begin{equation}  \label{eq:pauliJ}
J = \begin{bmatrix}
 0 & \qi \\ - \qi & 0  
 \end{bmatrix}.
 \end{equation}
\end{proposition-and-definition}

\medskip

The proof of Proposition~\ref{prop:11} is given in Section~\ref{section:3} below.  
As we will see there, the formula used to define $H$ in (\ref{eqn:H_Definition}) 
is based on the explicit rational formula for the $M_2 (\bC)$-valued 
$R$-transform of $Y$.

\begin{remark}
$1^o$ One can use \eqref{eqn:H_Definition} to interpret
our subordination function $\Omega$ as a subordination function for the second 
convolution power $\mu \boxplus \mu$ of some matrix valued distribution $\mu$.
Convolution powers of operator-valued distributions have been studied in 
\cite{AnshelevichBelinschiFevrierNica}. Since the focus of this paper is on explicit
computations of the Brown measure, we do not pursue this more abstract connection further here.

\vspace{6pt}

\noindent
$2^o$ For readers who are familiar with  the PDE method for Brown measures, below we list a dictionary 
of items which appear throughout the paper together with their translations into the
notation used in the PDE method.
\begin{align*}   P_\delta&=\varphi\left(\frac{1}{|\lambda-x|^2+\delta^2}\right),
   \\ P_\lambda&=\varphi\left(\frac{\overline{\lambda}-x}{|\lambda-x|^2+\delta^2}\right),
   \\ P_{\overline{\lambda}}&=\varphi\left(\frac{\lambda-x}{|\lambda-x|^2+\delta^2}\right).
\end{align*}
Moreover with the above notation we have
\begin{align*}
G_X( \, B(\lambda,\qi \delta) \, ) :=\begin{bmatrix}
        -\qi \delta P_\delta&P_{\overline{\lambda}}\\
        P_\lambda& -\qi \delta P_\delta
    \end{bmatrix}
\end{align*}
and
\begin{align}
    H_{11}(\lambda,\qi \delta)&:=H_{11}(B(\lambda,\qi \delta))=\qi \delta \left(1-p\frac{P_\delta}{|P_\lambda+i|^2+\delta^2P_\delta^2}
    \right), \label{eqn:H11-formula}\\
    H_{12}(\lambda,\qi \delta)&:=H_{12}(B(\lambda,\qi \delta))=\lambda-p\frac{P_{\overline{\lambda}}-\qi}{|P_\lambda+i|^2+\delta^2P_\delta^2}. \nonumber
\end{align}
\end{remark}

\bigskip

\subsection{Definition of the reparametrization $h : \cD \to \cM$.}
\label{subsection:12}
\phantom{x}

\noindent
Towards describing the reparametrization $h$ announced earlier in the introduction
we use matrices in $\bH_{2}^{+}$ which are of the special form
$B(\lambda, \zeta)$, with $\lambda \in \bC$ and $\zeta \in \bC^+$, as introduced
in Equation (\ref{B-Matrix}) above.
For such a matrix $B(\lambda, \zeta)$ and for $1 \leq i,j \leq 2$
we use the shorthand notation $H_{ij} ( \lambda , \zeta )$ 
for the $(i,j)$-entry of
$H \bigl(  B( \lambda , \zeta )  \bigr) \in M_2 ( \bC )$. 

We next record some special properties of $H_{11} ( \lambda, \zeta )$. 
\medskip

\begin{proposition}    \label{prop:12}
Fix $\lambda\in\mathbb{C}$.  Then:
\begin{enumerate}[$1^o$]
\item There exists a probability measure $\rho$ on $\mathbb{R}$, depending on $\lambda$, 
such that $H_{11}$ has an integral representation
\[
H_{11} ( \lambda, \zeta ) = \zeta + p\cdot\int_\mathbb{R} \frac{1}{\zeta-t} d \rho (t),
\ \ \forall \, \zeta \in \bC^{+}. 
\]
\item  $\Re H_{11} ( \lambda, \qi \delta ) = 0$,
$\forall \, \delta \in (0, \infty )$.
\item  The map 
\begin{equation}   \label{eqn:12a}
\delta \mapsto \frac {\Im  H_{11} ( \lambda,\qi \delta ) }{ \delta}
\end{equation}
is continuous and strictly increasing on $(0,\infty)$, and
has $\lim_{\delta \to \infty} 
\tfrac{ \Im  H_{11}( \lambda, \qi \delta ) }{\delta} = 1$.
\end{enumerate}
\end{proposition}

\begin{remark}
 Note that  $H_{22}=H_{11}$ and $H_{21}=\overline{H_{12}}$ and hence it follows
    from the preceding proposition that for fixed $\varepsilon>0$ and $\lambda\in\bC$,
    there is a unique $\delta>0$ such that the imaginary part of the matrix $H\left(\begin{bmatrix}
        \qi \delta&\lambda\\\overline{\lambda}& \qi \delta
    \end{bmatrix}\right)$ equals $\varepsilon I$. Since $\begin{bmatrix}
        \qi \delta&\lambda\\\overline{\lambda}& \qi \delta
    \end{bmatrix}$ is a point in the image of $\Omega$, this can be restated by saying that for fixed $\varepsilon>0$ the imaginary part of the matrix $\Omega\left(\begin{bmatrix}
        \qi \varepsilon & z\\\overline{z}& \qi \varepsilon
    \end{bmatrix}\right)$ is a function of the real part of this matrix. Such geometric properties of the subordination function  were studied in \cite{BelinschiGeometric}  in a much more general framework.
\end{remark}
\medskip

Proposition~\ref{prop:12} will be restated and proved in Section~\ref{sec:3.1} 
below.  Right now let us accept this proposition and observe that, as implied by its 
part $3^o$, the map (\ref{eqn:12a}) is sure to also have a limit when $\delta \to 0$.
We use this in order to define the desired pieces of structure $\cD, h, \cM$
appearing in our reparametrization, as follows.

\medskip

\begin{definition-and-remark}   \label{def:13}
Consider the same framework as above.  
\begin{enumerate}[$1^o$]
    \item 
 We let
\begin{equation}   \label{eqn:13a}
\cD:=\left\{  \lambda\in\bC \mid  \lim_{\delta\to 0^+}
\frac{ \Im  H_{11} ( \lambda, \qi \delta ) }{\delta} < 0 \right\} .
\end{equation}
Due to $3^o$ of Proposition~\ref{prop:12}, $\cD$ is an open subset of $\bC$.

\item 
Fix $\lambda\in\cD$. The properties of the function indicated in (\ref{eqn:12a})
ensure that there exists a number $\delta_0 ( \lambda ) \in (0,\infty)$, uniquely determined,
such that $\tfrac{ \Im  H_{11} ( \lambda, \qi \delta_0 (\lambda)) }{\delta_0(\lambda)} = 0$. 
In view of Proposition $\ref{prop:12}.2^{o}$, this number can be equivalently described by 
its property that
\begin{equation}   \label{eqn:13b}
H_{11} \bigl( \, \lambda, \qi \delta_0(\lambda) \, \bigr) = 0.
\end{equation}
We define
\begin{equation}   \label{eqn:13c}
h( \lambda ) := H_{12} \bigl( \, \lambda, \qi \delta_0 (\lambda) \, \bigr).
\end{equation}

\item We put $\cM := h ( \cD ) \subseteq \bC$.
\end{enumerate}
\end{definition-and-remark}

\medskip

\begin{remark}   \label{rem:14}
\begin{enumerate}[$1^o$]
    \item []
    \item
 Part $2^o$ of the preceding definition introduces, at the same time with $h$,
a function $\delta_0 : \cD \to (0, \infty)$.  The function $\delta_0$ will play a significant 
role in what follows.  We record here that the relations defining $\delta_0$ and $h$ can be 
consolidated in the formula
\begin{align}    \label{eqn:H_delta_0}
H\left( \, \begin{bmatrix} \qi \, \delta_0 (\lambda) & \lambda\\
		\overline{\lambda} & \qi \, \delta_0 (\lambda) \end{bmatrix} \, \right) 
=\begin{bmatrix} 0 & h(\lambda)\\
		\overline{h(\lambda)} & 0\end{bmatrix}, \ \ \lambda \in \cD.
\end{align}
\item 
In order to fix the ideas of how $\cD$ and $\cM$ may look like, Figure 1
shows them in the case when $\pp =1$ and $x$ has symmetric Bernoulli distribution.
This example is worked out in full detail in Section~\ref{subsection:14} below.
\end{enumerate}

\begin{figure}[ht]
	\centering
	\begin{subfigure}[c]{0.4\textwidth}
		\centering
		\includegraphics[width=\textwidth]{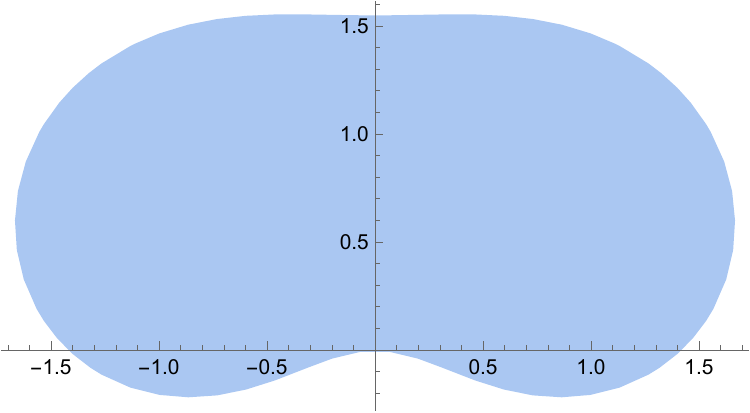}
		\caption{The set $\cD$.}
		\label{fig:sub1}\end{subfigure}
	\begin{subfigure}[c]{0.25\textwidth}
		\centering
		\includegraphics[width=\textwidth]{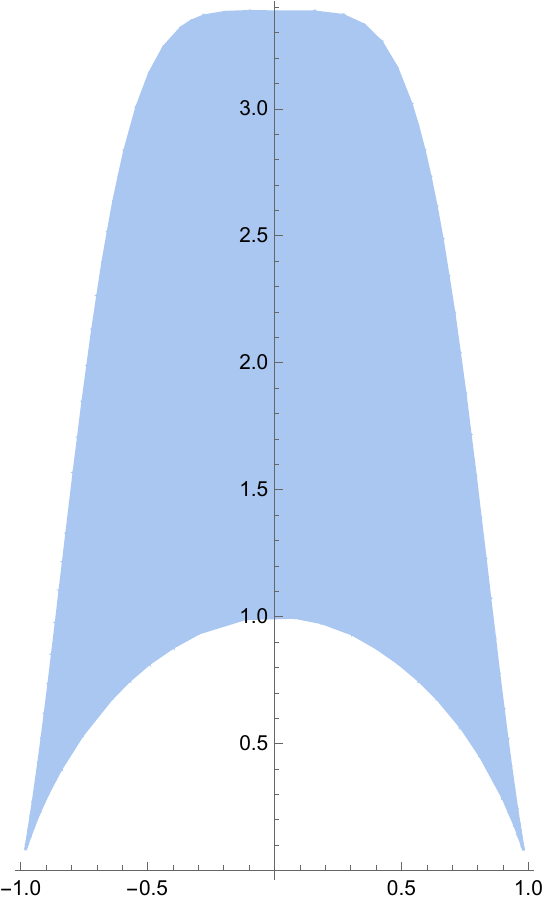}
		\caption{The set $\cM$.}
		\label{fig:sub2}
	\end{subfigure}
	    \caption{The sets $\cD$ and $\cM$ 
        in the case when $\pp = 1$ and $x$ has symmetric Bernoulli distribution. One gets 
        $\cD = \bigl\{ \lambda \in \bC : 1+| \lambda |^2 > | 1 - \lambda ( \lambda - \qi ) |^2 \, \bigr\}$,
        or, in cartesian coordinates:  \\
        $\cD = \bigl\{ (\alpha,\beta)\in\bR^2 : \left(\alpha ^2+(1-\beta ) \beta +1\right)^2+(\alpha  (\beta -1)+\alpha  \beta )^2-\alpha ^2-\beta ^2-1<0 \, \bigr\}$.
        Then
       $\cM=\{(s,t)\in \bR^2:4 \left(s^2-1\right) t^3+2 \left(s^2+1\right) \left(s^2-1\right)^2 t^2+\left(s^2-1\right)^4+\left(s^4-s^2+1\right) t^4<0\}$} 
        
    \label{fig:combined}
\end{figure}
\end{remark}

\bigskip

\subsection{Properties of $h : \cD \to \cM$, and the main theorem.}
\label{subsection:13}
\phantom{x}

\begin{proposition}    \label{prop:15}
\begin{enumerate}[$1^o$]
    \item []
    \item  $\cD$ and $\cM$ are open subsets of $\bC$.
    \item $\delta_0$ and $h$ are differentiable on $\cD$.
    \item  $h$ is injective (hence a bijection from $\cD$ onto $\cM$).
    \item The function $h^{-1}$ is differentiable on $\cM$.
\end{enumerate}
\end{proposition}

The properties of $\delta_0$ and $h$ gathered in Proposition \ref{prop:15} 
are proved at various points of the Sections 3 and 5 of the paper. 

The main result of the paper will be stated under the following assumption, which in particular rules 
out the delicate situation when the boundary of $\cD$ contains atoms of $\mu_x$. 

\medskip

\begin{assumption}
  \label{assumption:16}
Recall (cf.~the description of framework on page 3) that $\mu_x$ denotes the distribution 
of $x$ with respect to $\varphi$.  We assume that $\mu_x$ is not a point mass, that the support of the 
continuous part of $\mu_x$ is contained in $\cD$, and that there are no atoms of $\mu_x$ on the boundary 
$\mathrm{bd} ( \cD ) = \mathrm{cl}(\cD)\setminus \cD$.
\end{assumption} 

\medskip

\begin{remark}
In Proposition \ref{prop:43} below we provide some sufficient conditions on 
the Cauchy transform of $\mu_x$ which ensure
that Assumption~\ref{assumption:16} holds.
\end{remark}

\vspace{6pt}

\begin{remark}    \label{remark:1.9}
We note that for $\lambda\notin \spec(x)$ we have 
\begin{align*}
     \lim_{\delta\rightarrow {0}}\det(J-G_X(B(\lambda,\qi \delta))&=\det(J-G_X(B(\lambda,0))\\
      &=\begin{vmatrix}
          0 & \qi - \overline{G_{\mu_x}(\lambda)}\\
          -\qi-G_{\mu_x}(\lambda) & 0
      \end{vmatrix}=-|\qi+G_{\mu_x}(\lambda)|^2,
\end{align*}
and \eqref{eqn:H_Definition} tells that 
\begin{equation}
   \label{eqn:limit-H11-over-delta}
      \lim_{\delta\rightarrow{0^+}}\frac{\Im H_{11}(B(\lambda,\qi\delta))}{\delta}
      =1-\frac{p\cdot \varphi \bigl( \,  |\lambda-x|^{-2} \, \bigr) }{ |i+G_{\mu_x}(\lambda)|^2}.
\end{equation}
By the definition of $\cD$, for $\lambda\not\in\cD\cup\spec(x)$ we have that the right-hand side of 
(\ref{eqn:limit-H11-over-delta}) is $\geq 0$, which implies
\begin{equation}
  \label{eqn:inequality-outsidep}
  \varphi\left( \,|\lambda-x|^{-2} \, \right)
  \cdot \frac{1}{|\qi+ G_{\mu_x}(\lambda)|^2}\leq \frac{1}{p}. 
\end{equation}
\\Define 
$F(\lambda)= |i+G_{\mu_x}(\lambda)|^2 \left( \, \varphi\left( |\lambda-x|^{-2} \right) \, \right)^{-1}$ for 
$\lambda\not\in\spec(x)$.
In Section 4, we show that $F$ can be extended continuously to $\spec(x)$ and then set $\cD$ can be described simply by \[\cD=\{\lambda\in\bC:F(\lambda)<p\}.\]
Similarly one can check that for $\lambda\not\in\spec(\cD)$
\[\lim_{\delta\to 0^+} H_{12}(\lambda,\qi \delta)=\lambda-p\frac{1}{\qi+G_{\mu_x}(\lambda)}.\]
Later in the paper (see Subsection \ref{subsec:extension}) we extend $\delta_0$ to $\bC$ by $\delta_0(\lambda)=0$ for $\lambda\not\in \cD$ and then the right-hand side of the equality above gives similar extension of $h$ outside $\cD$. Among other results we show that such $h$ is a homeomorphism of the complex plane. Note however that this mapping is holomorphic only outside $\mathrm{cl}(\cD)$.
\end{remark}

\noindent
\medskip

\begin{theorem}   \label{thm:18}  Consider Assumption~\ref{assumption:16}
and the setting described above (where, in particular, the map  $h$ from Definition~\ref{def:13}
has an extension to a homeomorphism $h : \bC \to \bC$). Then the map $h$ and the Brown measure $\mu_{x+\qi y}$ have the following properties.

$1^o$ Every atom $\alpha\in\spec(x)\setminus \mathrm{cl}(\cD)$ has mass $\mu_x(\{\alpha\})>p$ and is fixed by $h$, i.e.,  $h(\alpha)=\alpha$. 

$2^o$ Let $\alpha_1,\ldots,\alpha_k$ be the set of atoms of $\mu_x$ in $\spec(x)\setminus \cD$, then the support of $\mu_{x+\qi y}$ is contained in $\mathrm{cl}(\cM)\cup\{\alpha_1,\ldots,\alpha_k\}$. 

$3^o$ For any $i\in\{1,\ldots,k\}$ we have $\mu_{x+\qi y}(\{\alpha_i\})=\mu_x(\{\alpha_i\})-p$.

$4^o$  On $\cM$, $\mu_{x+\qi y}$ has a density $f$  which for any $s+\qi t\in\cM$ is given by the formula 
\begin{equation}   \label{eqn:density}
        f(s,t)=\frac{1}{2 \pi}\left[\frac{1}{t}\left(\frac{\partial \alpha}{\partial s}
        +\frac{\partial \beta}{\partial t}\right)-\frac{1}{t}-\frac{\beta}{t^2}\right],
\end{equation}
where $\alpha+\qi\beta=h^{-1}(s+\qi t)$.

\end{theorem}

\begin{remark}
    The formula for the Brown measure stated above can also be understood in terms of the derivative of the logarithm of the Fuglede-Kadison determinant (which we introduce in Section 2 below). From this point of view our result is similar to several other recent results on Brown Measure (see \cite{BelinschiYinZhong2024, BercoviciZhong2022, DHK2022Brown, HallHo2022LMP, HoZhong2023Brown, Zhong2021Brown}), namely for every $z\in\bC$ we have 
    \[\frac{\partial L_{x+\qi y}}{\partial z}=\varphi\left(\frac{(\lambda-x)^*}{|\lambda-x|^2+\delta_0(\lambda)^2}\right), \quad\mbox{with } \lambda=h^{-1}(z),\]
    where $\delta_0(\lambda)$ is regarded as subordination function in free probability or the canonical choice of $\delta$ in the Hamilton-Jacobi method for Brown measure. 
\end{remark}

\begin{remark}
    In the physics literature, Jarosz and Nowak \cite{JaroszNowak} give a general algorithm 
    for computing Brown measures of operators of the form $x+\qi y$ for $x,y$ selfadjoint
    and free. While the method of Jarosz and Nowak is not mathematically rigorous as it is
    written, and the example of $y$ free Poisson is not explicitly treated in \cite{JaroszNowak}, 
    we expect that a detailed implementation of their algorithm in this setting might lead to a 
    result which is equivalent to the one obtained above. This is analogous to how things went 
    in the case when $y$ is semicircular, where the agreement between \cite{JaroszNowak} and rigorous 
    analysis was demonstrated in detail in \cite[Section 9]{HallHo2022LMP}. The main difficulty 
    in relating our approach to the algorithm of Jarosz and Nowak is that the method
    used in \cite{JaroszNowak} does not involve the reparametrization mapping $h$. 
    For a special case verification, we mention that in the case when $x$ has symmetric Bernoulli 
    distribution and $y$ is free Poisson, one can easily verify that the algorithm from 
    \cite{JaroszNowak} yields the same explicit density as the one we state in the next 
    Subsection \ref{subsection:14}.
\end{remark}

\begin{remark}
    Observe that in Theorem \ref{thm:18} we do not describe the Brown measure of the boundary of $\cM$, 
    however we conjecture that the Brown measure of the boundary of $\cM$ is zero. Also we do not show 
    that the density of the Brown measure is strictly positive at every point of $\cM$.
\end{remark}

\bigskip

\subsection{A worked example: $x$ symmetric Bernoulli, $y$ standard free Poisson.}
\label{subsection:14}

\phantom{x}

\noindent
As announced in Remark~\ref{rem:14} above, we now present what the 
methodology described in the preceding subsections will concretely amount to, in the
special case when the element $x$ has symmetric Bernoulli distribution 
$\mu_x = \frac{1}{2} \bigl( \delta_1 + \delta_{-1} \bigr)$, and the parameter $\pp$ 
of the free Poisson element $y$ is equal to $1$.
More precisely, we carry out the following computations:
\begin{enumerate}[$1^o$]
    \item Explicit formulas of the left inverse $H(B)$;
    \item An algebraic description of the domain $\cD$;
    \item Explicit formulas of the functions $h$ and $\delta_0$;
    \item An implicit algebraic equation for the support of the Brown measure;
    \item An implicit algebraic equation for the spectrum, which is thus shown to coincide with the support of the Brown measure.
\end{enumerate}

{$1^o$}
Due to the special choice for the distribution of $x$, the Cauchy transform of the matrix 
$X \in M_2 ( \cA )$ from Equation (\ref{eqn:1b}) comes to
\[
G_X (B) = \frac{1}{2} \Bigl( \, (B-K)^{-1} + (B+K)^{-1} \, \Bigr), \ \ B \in \bH_2^{+},
\]
where we denoted $K = \begin{bmatrix} 0 & 1 \\ 1 & 0 \end{bmatrix}$. 
It follows from $G_X (B) = \int_{\bR} (B - tK)^{-1} \, d \mu_x (t)$, (see \cite[p.29]{BelinschiMaiSpeicher}).
The function $H$ introduced in Definition~\ref{prop:11} will thus take here the form 
\begin{equation*}
H(B) = B+\left(J-\frac{1}{2} 
\left( \left( B-K \right)^{-1} + \left( B+K \right)^{-1} \right) \right)^{-1},
\ \ B \in \bH_2^{+}.
\end{equation*}

{$2^o$}
For $B=\begin{bmatrix}
    \qi \delta&\lambda\\
    \overline{\lambda}&  \qi \delta  \
\end{bmatrix}$, the explicit formula for the $(1,1)$-entry $H_{11} (\lambda , \qi \delta)$ 
of $H(B)$ comes out as follows:
\begin{equation*}
H_{11} (\lambda , \qi \delta ) = \qi \delta \left(1 -\frac{  
\delta ^2 + \lambda  
\overline{\lambda}+1
}{\delta ^4+\delta ^2 
(\lambda  (2 \overline{\lambda}+\qi)-\qi \overline{\lambda}+3)
+(-1+\lambda  (\lambda -\qi)) (-1+\overline{\lambda} (\overline{\lambda}+\qi))}\right)
\end{equation*}

Thus, we have 
\begin{equation}   \label{eqn:1-4x}
\frac{\Im H_{11} (\lambda, \qi \delta) }{\delta} = 1 -
\frac{ \delta ^2+\lambda  \overline{\lambda}+1 }{\delta ^4+\delta ^2 (\lambda  (2 \overline{\lambda}+\qi)
-\qi \overline{\lambda}+3)+(-1+\lambda  (\lambda -\qi)) (-1+\overline{\lambda} (\overline{\lambda}+\qi))}.
\end{equation}
Going next to the definition of $\cD$ in Equation (\ref{eqn:13a}), we see that the $\delta \to 0^+$
limit considered there is found, in this particular case, by simply substituting $\delta = 0$ on the
right-hand side of (\ref{eqn:1-4x}).   This gives
\[
\cD = \{\lambda\in\bC  :  \ 1 - \frac{ |\lambda|^2 + 1 }{ (1-\lambda  (\lambda -\qi)) 
                                    (1-\overline{\lambda} (\overline{\lambda}+\qi)) } <0 \},
\]
which is precisely the domain shown in Figure 1(A).

{$3^o$} Our next goal is to determine the function $h$. First we have to determine $\delta_0(\lambda)$, 
that is, we need to find the solution to the equation 
\[
1 -\frac{  
\delta ^2+\lambda  \overline{\lambda}+1}{\delta ^4
          + \delta ^2 (\lambda  (2 \overline{\lambda}+\qi)-\qi \overline{\lambda}+3)
          + (-1+\lambda  (\lambda -\qi)) (-1+\overline{\lambda} (\overline{\lambda}+\qi)) } = 0.
\]
The latter equation can be re-written as an algebraic equation of degree 4 in $\delta$.
Among its four solutions, only one maps the imaginary axis to the real axis (which is a necessary 
property of the mapping $\delta_0$) and we obtain 
\[
\delta_0(\lambda) = \frac{\sqrt{\sqrt{3 \lambda ^2+10 \lambda  \overline{\lambda}
+3 \overline{\lambda}^2+4}-2 \lambda  \overline{\lambda}-\qi \lambda +\qi \overline{\lambda}-2}}{\sqrt{2}}.
\]

Since $h(\lambda)$ is the $(1,2)$-entry in $H \left(\begin{bmatrix}
    \qi\delta_0(\lambda)&\lambda\\
    \overline{\lambda}&\qi\delta_0(\lambda)
\end{bmatrix}\right)$, substituting the above $\delta_0 ( \lambda )$ into the explicit 
formula for $H$ gives
\begin{align*}
    h(\lambda)=\frac{\sqrt{(3 \lambda +\overline{\lambda}) (\lambda +3 \overline{\lambda})+4}
    +\lambda  (2 \lambda -\qi)+\overline{\lambda} (-2 \overline{\lambda}+\qi)}{2 (\lambda +\overline{\lambda}-\qi)}.
\end{align*}

{$4^o$}
Based on the above we find the equation of the boundary of $\cM=h(\cD)$ which comes out as shown in the caption of Figure \ref{fig:combined}. We repeat Figure \ref{fig:combined}(b), by also showing the eigenvalues of a sample of a random matrix model, confirming the
heuristically expected result that Brown measure is approximated by eigenvalues of non-normal
matrices. Here we took matrices of size $N\times N$ with $N=2000$ and $X_N$ diagonal with half
eigenvalues $-1$ and half $1$ and $Y_N$ being the Wishart random matrix that is 
$Y_N=\frac{1}{N}G_N G_N^*$ where $G_N$ is a GUE matrix and calculate the eigenvalues of 
$X_N+\qi Y_N$.

\begin{figure}[ht]
	\centering
	\begin{subfigure}[c]{0.3\textwidth}
		\centering
		\includegraphics[width=\textwidth]{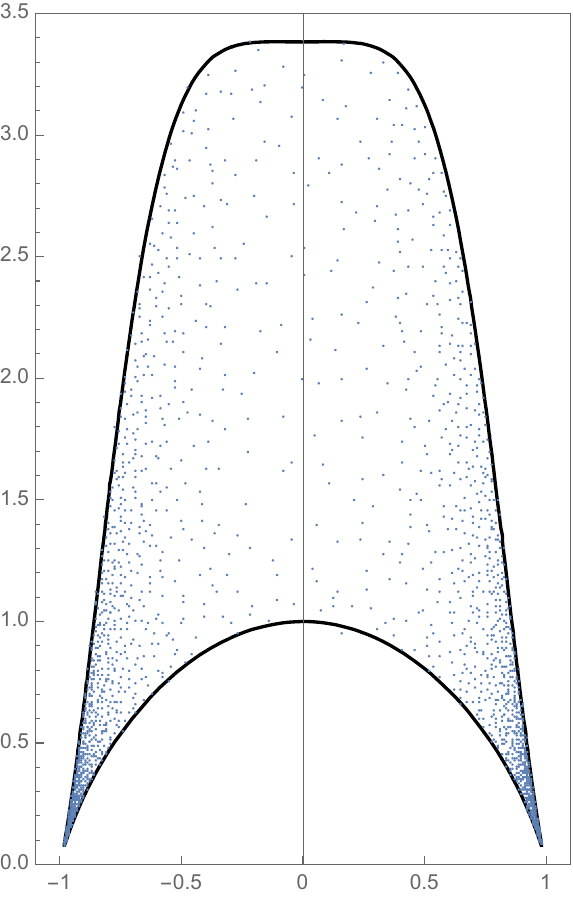}
		\caption{Eigenvalue simulation}
        \end{subfigure}
	\begin{subfigure}[c]{0.315\textwidth}
		\centering
		\includegraphics[width=\textwidth]{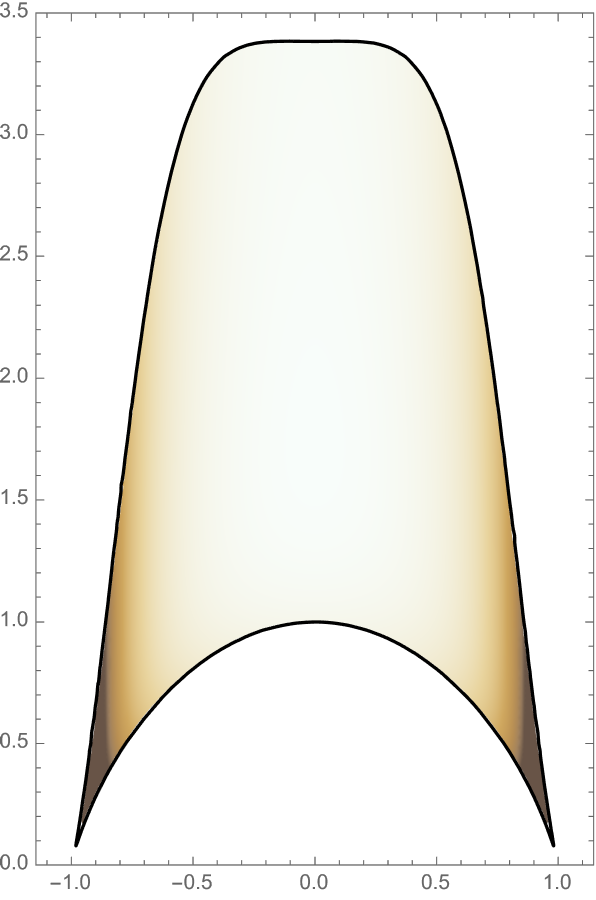}
		\caption{Density plot}
	\end{subfigure}
    \begin{subfigure}[c]{0.25\textwidth}
		\centering
		\includegraphics[width=\textwidth]{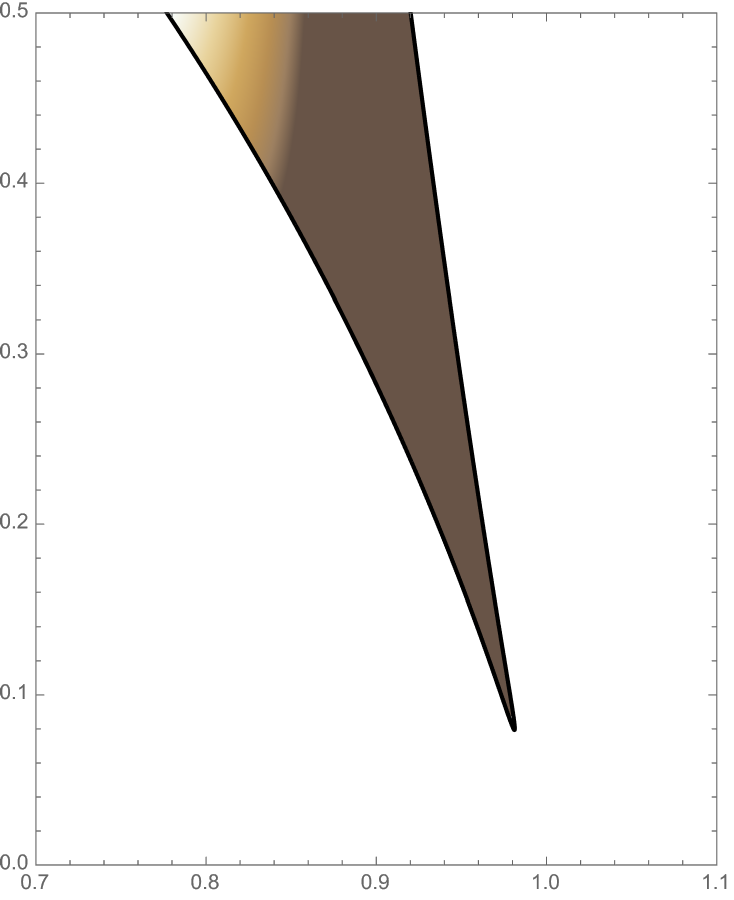}
		\caption{Density near $(1,0)$}
	\end{subfigure}
	    \caption{The black line represents the boundary of $\cM$.
    The blue dots in (A) are the eigenvalues of the matrix approximation $X_N+\qi Y_n$.} 
        
    \label{fig:combined}
\end{figure}

To get the explicit formula for the density one has to invert $h(\alpha+\qi \beta)=s+\qi t$,
that is write $\alpha$ and $\beta$ as functions of $s$ and $t$.
We see that for $\alpha+\qi \beta\in\cD$ we have 
\[
h(\alpha+\qi\beta)
= \frac{ \sqrt{4 \alpha ^2+\beta ^2+1}+4 \qi \alpha  \beta +\beta }{ 2 \alpha - \qi }. 
\]
Separating real and imaginary part we arrive at the system of equations
\begin{equation*}
\left\{
    \begin{aligned}
        s&=\frac{2 \alpha  \left(\sqrt{4 \alpha ^2+\beta ^2+1}-\beta \right)}{4 \alpha ^2+1}\\
        t&=\frac{\sqrt{4 \alpha ^2+\beta ^2+1}+8 \alpha ^2 \beta +\beta }{4 \alpha ^2+1}.
    \end{aligned}
    \right.
\end{equation*}
Solving this system we obtain
\begin{align}\label{eqn:alphabeta}
\alpha = \frac{s t}{2 \left(1-s^2\right)}, \ \mbox{ and }
\ \beta = \frac{s^2+t^2-1}{2 t}.
\end{align}

From Theorem~\ref{thm:18} we obtain that for $s+\qi t\in \cM$ the density 
of the Brown measure of the element $x+\qi y$ is given by

\[
f(s,t) = \frac{1}{2 \pi} \cdot \frac{1}{t} 
\left( \frac{1 - s^2}{t^2} + \frac{(1 + s^2) t}{2 \left( 1 - s^2 \right)^2} - 1 \right) .
\]
While this calculation is carried in the case $p=1$ we point out that one can allow general parameter $p>0$ for the free Poisson element $y$. Then equations \eqref{eqn:alphabeta} take the form
\begin{align*}
\alpha = \frac{s t}{2 \left(1-s^2\right)}, \ \mbox{ and }
\ \beta = \frac{s^2+t^2+(1-p)t-1}{2 t},
\end{align*}
and thus one gets again an explicit density on appropriate set $\cM$, moreover for $p<1/2$ points $\pm 1$ are outside $\cD$ and the Brown measure of $x+\qi y$ has atoms at $\pm 1$ both of mass $1/2-p$.

{$5^o$} We conclude by showing that
   in the present example the support of the Brown measure coincides with the spectrum. 
Indeed using resolvent estimates from \cite{Lehner:2001:computation}
one can compute the outer boundary of the spectrum of $x+\qi y$ and show
that it is equal to the boundary of the support of the Brown measure.
More precisely, given arbitrary (not necessarily selfadjoint) $*$-free
random variables $a_1$ and $a_2$,
denote by  $f(\zeta_1) = M_{a_1}(\zeta_1)$ and $g(\zeta_2)=M_{a_2}(\zeta_2)$ 
their moment generating functions.
Pick numbers $\zeta_1,\zeta_2$ in the respective domains of convergence such that $\zeta_1f(\zeta_1)=\zeta_2g(\zeta_2)\ne 0$.
Then  the resolvent of $a_1+a_2$ at 
\begin{equation}
\label{eq:lambda=1s+1t-1sf}
    z=\frac{1}{\zeta_1} + \frac{1}{\zeta_2} - \frac{1}{\zeta_1f(\zeta_1)}
\end{equation}
can be written as
\begin{equation}
      z-a_1-a_2
    = \frac{g(\zeta_2)}{\zeta_1}
      (1-\zeta_1a_1)
      \left(
        1 - \frac{\bub{a}_1(\zeta_1)\,
                  \bub{a}_2(\zeta_2)}
                 {f(\zeta_1)\, g(\zeta_2)}
      \right)
      (1-\zeta_2a_2)
\end{equation}
where $\bub{a}_1(\zeta_1)=(1-\zeta_1a_1)^{-1}-f(\zeta_1)$ and  $\bub{a}_2(\zeta_2)=(1-\zeta_2a_2)^{-1}-g(\zeta_2)$
are the centered resolvents of $a_1$ and $a_2$.

Now a result of Haagerup and Larsen \cite[Proposition~4.1]{HaagerupLarsen} asserts that the spectral radius of a product of two free centered random variables is equal to its $L^2$-norm and moreover the spectrum
includes the entire circle \cite[Corollary~4.1]{Lehner:2001:computation}.
Then \cite[Proposition~4.2]{Lehner:2001:computation} provides the
following criteria:
\begin{enumerate}[(i)]
\item 
  \label{eq:normXnormY<1}
  if $\norm{\bub{a}_1(\zeta_1)}_2 \norm{\bub{a}_2(\zeta_2)}_2 < \abs{f(\zeta_1)g(\zeta_2)}$
  then $z\not\in\spec(a_1+a_2)$
  \item 
    \label{eq:normXnormY=1}
 if $\norm{\bub{a}_1(\zeta_1)}_2 \norm{\bub{a}_2(\zeta_2)}_2 = \abs{f(\zeta_1)g(\zeta_2)}$
  then  $z\in\bd(\spec(a_1+a_2))$
  unless  \eqref{eq:lambda=1s+1t-1sf}
 has a solution satisfying \eqref{eq:normXnormY<1}.
\end{enumerate}

That is, the region
$$
\bigl\{  z = \frac{1}{\zeta_1} + \frac{1}{\zeta_2} - \frac{1}{\zeta_1f(\zeta_1)}
\mid   \norm{\bub{a}_1(\zeta_1)}_2 \norm{\bub{a}_2(\zeta_2)}_2 < \abs{f(\zeta_1)g(\zeta_2)}
\bigr\}
$$
is part of the resolvent set of $a_1+a_2$ and its boundary
is part of the curve
\begin{equation}
\label{eq:bdspec}
\bigl\{  z = \frac{1}{\zeta_1} + \frac{1}{\zeta_2} - \frac{1}{\zeta_1f(\zeta_1)}
\mid   \norm{\bub{a}_1(\zeta_1)}_2 \norm{\bub{a}_2(\zeta_2)}_2 = \abs{f(\zeta_1)g(\zeta_2)}
\bigr\}
\end{equation}
and contained in the spectrum of $a_1+a_2$.

Carrying out the calculations for $a_1=x$ and $a_2=\qi y$,
where $f(s)=\frac{1}{1-s^2}$ and $g(t)$ satisfies the quadratic equation
$itg(t)^2-g(t)+1=0$, one arrives after a few elimination steps at the following algebraic equation 
for the curve \eqref{eq:bdspec}

\begin{equation}
\label{eq:boundaryspec}
s^8+2 s^6 t^2+s^4 t ^4 -4 s^6 -2 s^4 t^2 -s ^2 t^4+4 s^2 t^3+6 s ^4 -2 s^2 t^2+t^4 -4 t ^3 -4 s^2+2 t^2+1 
=0
.
\end{equation}
This is an algebraic curve of genus 1 (in fact elliptic) and has two isolated points
at $\pm1$.
These however are not in the spectrum because for $\lambda=\pm1$ there are solutions $\zeta_1$ and $\zeta_2$
of \eqref{eq:lambda=1s+1t-1sf} such that 
   $\norm{\bub{a}_1(\zeta_1)}_2 \norm{\bub{a}_2(\zeta_2)}_2 < \abs{f(\zeta_1)g(\zeta_2)}$ and thus the points $\pm1$ belong to the resolvent set.
The rest of the curve however coincides with the boundary of the region $\cM$ from
Fig.~\ref{fig:combined}
and thus the support of the Brown measure fills the full spectrum.

\vspace{10pt}

\subsection{Comments around some technical aspects of the paper.}
\phantom{x}

\noindent
The explicit description of the density requires several technical steps not described above.

\smallskip
\noindent
\begin{itemize}
    \item [--]
 Taking advantage of the explicit form of the left inverse $H$, 
we extend $\Omega$ to a neighborhood of 
$B=\begin{bmatrix} 0&z\\ \overline{z}&0 \end{bmatrix}$ for every $z\in \cM$. This step relies on the Inverse Function Theorem, which in turn means that we need to prove that some Jacobians are nonzero.
\item[--] 
The above implies that we can extend the Cauchy transform of $X+Y$ and for $z\in \cM$ have
\[G_{X+Y}\left(\begin{bmatrix}
    0&z\\\overline{z}&0
\end{bmatrix}\right)=G_X\left(\begin{bmatrix}
    \qi \delta_0(\lambda)&\lambda\\\overline{\lambda}&\qi \delta_0(\lambda).
\end{bmatrix}\right)\]
with $z=h(\lambda)$ and $\delta_0(\lambda)>0$.
\item[--]
 Relations between $L_{x+iy}$ and the operator valued Cauchy transform, and the above regularity of $G_{X+Y}$ allow us to show that the function 
\[ (s,t,\varepsilon)\mapsto L_{x+\qi y}( s+\qi t,\varepsilon)
 \]
extends to a real analytic function in some neighborhood of $(s, t, 0)$ for $s+\qi t\in\cM$.
\item[--]
 In order to conclude that the Brown measure is supported inside $\mathrm{cl}(\cM)$ we extend $\delta_0$ and $h$ beyond $\mathrm{cl}(\cD)$ and show that $h$ is holomorphic there.

 \end{itemize}
 
\subsection{Organization of the paper}

\phantom{x}

\noindent
Besides this introduction, the paper has 5 other sections:

-- In Section~2 we present background. 

-- Section~3 provides the details of the function $H_{11}$ and properties of $\delta_0$ 
and $h$ inside and outside $\cD$. 

-- In Section 4 we discuss Assumption \ref{assumption:16}.

-- In Section~5 we study positivity of the Jacobians inside and outside the domain $\cD$. 
s
-- In Section 6 we show real analyticity of $L_{x+\qi y}$, study the support of Brown measure and derive the formula 
for the density of the Brown measure.

-- In Section 7 we study discrete part of the Brown measure and prove Theorem \ref{thm:18}

\section{Background}  \label{section:2}

\subsection{Brown measure}
\phantom{x}

\noindent
Let $(\mathcal{A},\varphi)$ be a $W^*$-probability space, and let $a\in \mathcal{A}$. The Fuglede-Kadison 
determinant $\KFD(a)$ is defined by 
\[
  \log \KFD(a)=\varphi(\log(|a|)),
\]
where $|a|=\sqrt{a^*a}$ by functional calculus.  
Brown showed in \cite{Brown} that the function $z\mapsto \log \KFD(z-a)$ is a subharmonic function on
the complex plane $\mathbb{C}$. The Brown measure of $a$ is the unique probability measure $\mu_a$ on $\mathbb{C}$ determined by the identity
\[
  \int_\mathbb{C}\log\vert z-\lambda \vert d\mu_a(\lambda)
   =\log \KFD(z-a).
\]
We denote the logarithmic moment by
\[
  L_a(z)=\varphi(\log(|z-a|)),
\]
and its regularization by
\[
 L_a(z,\varepsilon)=\frac{1}{2}\varphi(\log(|z-a|^2+\varepsilon^2)).
\]
The Brown measure can be calculated, in distributional sense, by 
\[
\mu_a= \frac{2}{\pi} \frac{\partial}{\partial z}\frac{\partial}{\partial \overline{z}} L_a(z).
\]
The regularized Brown measure of $a$ with parameter $\varepsilon$
is calculated as 
\[ 
\mu_{a,\varepsilon} =
\frac{2}{\pi} \frac{\partial}{\partial z} \frac{\partial}{\partial \overline{z}} L_a(z,\varepsilon).
\]
It is known that $\mu_{a,\varepsilon}\rightarrow{\mu_a}$ weakly as $\varepsilon\rightarrow{0}$. 
The reader is referred to \cite{Brown, HaagerupSchultz2007, BelinschiSniadySpeicher} for details. 

\vspace{10pt}

\subsection{Matrix valued subordination}\label{subsec:2.2}
\phantom{x}

\noindent
We will rely on results from \cite{Voiculescu2000, BelinschiMaiSpeicher} concerning the existence and properties of subordination functions in the operator-valued setting (cf.\ Theorem~2.2 of \cite{BelinschiMaiSpeicher}). We apply these results in the framework described in Subsection~\ref{subsection:11}; recall that \(X\) and \(Y\) are free with amalgamation over the algebra \(M_2(\mathbb{C})\).

In this framework, there exists a Fr\'echet-analytic map
\[
\Omega \colon \mathbb{H}_2^+ \to \mathbb{H}_2^+
\]
such that:
\begin{enumerate}
    \item \(\Im \Omega(B) \geq \Im B\) for all \(B \in \mathbb{H}_2^+\);
    \item \(G_{X+Y}(B) = G_X(\Omega(B))\) for all \(B \in \mathbb{H}_2^+\).
\end{enumerate}

In the next subsection, we study the subordination function \(\Omega\) and its 
left inverse in greater detail, focusing on the special case where \( y \) is a 
free Poisson random variable with parameter \(p > 0\).

\subsection{Free Poisson elements -- the subordination function and its left inverse}

\begin{lemma}  \label{lemma:RtransformY}
Let $y$ be a free Poisson random variable with rate $\theta=1$ and jump size $p>0$.
\begin{enumerate}[$1^o$]
    \item 
 The operator-valued free cumulants of $Y$ are given by
\begin{equation}   \label{eqn:31z}
\kappa_n^\cB (YB, \cdots, YB, Y)=p ( J B)^{n-1} J ,
\ \ \forall \, n \in \bN \mbox{ and } B\in M_2(\mathbb{C})
\end{equation}
where $J$ is the Pauli matrix \eqref{eq:pauliJ}.
\item
The $R$-transform of $Y$ is rational
$R_Y (B)=p (J- B)^{-1}$
and 
has an immediate analytic continuation to 
$\bH_2^{+} \cup  \bH_2^{-}$.
\end{enumerate}

\end{lemma}

\begin{proof}
$1^o$ We write $Y = Y_o J$ with
$Y_o := \begin{bmatrix} y &  0 \\ 0  & y \end{bmatrix}$, and we note that $Y_o$
commutes with every matrix in $M_2 ( \bC )$.  So then
$\kappa_n^\cB (YB, \cdots, YB, Y)
= \kappa_n^\cB \bigl(Y_o (JB), \cdots, Y_o (JB), Y_o J \bigr)$,
where we observe that on the right-hand side all 
the occurrences of the matrix $JB \in M_2 ( \bC)$ can be pulled towards 
left until they become pre-factors. For example,
in order to pull out the first two occurrences of $JB$ we go as follows:
\begin{align*}
\kappa_n^\cB \bigl( \, Y_o (JB), Y_o (JB), \cdots \bigr)
& = \kappa_n^{\cB} \bigl( \, (JB) Y_o, (JB) \, Y_o, \cdots \bigr)       \\
& = (JB) \cdot \kappa_n^{\cB} \bigl( \, Y_o, (JB) Y_o, \cdots \bigr)  
    \  \mbox{ (left module property of $\kappa_n^{\cB}$) }  \\
& = (JB) \cdot \kappa_n^{\cB} \bigl( \, Y_o \, (JB), \,  Y_o, \cdots )  
    \  \mbox{ (balancing property of $\kappa_n^{\cB}$) }     \\
& = (JB) \cdot \kappa_n^{\cB} ( \, (JB) \, Y_o , \, Y_o, \cdots )   \\ 
&  = (JB)^2 \cdot \kappa_n^{\cB} ( \, Y_o ,  Y_o, \cdots ). 
\end{align*}
The outcome of pulling all the $JB$ factors to the left is that
\begin{equation}   \label{eqn:31a}
\kappa_n^\cB (YB, \cdots, YB, Y)
= (JB)^{n-1} \cdot \kappa_n^\cB (Y_o, \cdots, Y_o, Y_o ) \cdot J.
\end{equation}
On the other hand, the special form of $Y_o$ immediately implies that
\[
\kappa_n^\cB (Y_o, \cdots, Y_o, Y_o ) 
= \begin{bmatrix} \kappa_n (y) &  0 \\ 0  & \kappa_n (y) \end{bmatrix}
= \kappa_n (y) \, I_2 = p I_2, \ \ n \in \bN,
\]
and substituting this into (\ref{eqn:31a}) gives the required formula
(\ref{eqn:31z}).

\medskip

\noindent
$2^o$ For $\Vert B\Vert< 1$ we have 
\begin{align*}
R_Y(B) = \sum_{n\geq 0} \kappa_{n+1}^\cB(YB, \cdots, YB, Y)
  =\sum_{n\geq 0} p(JB)^{n}J & =p (J- B)^{-1}.
\end{align*}
Note that for $B$ such that $\Im B<0$ we have $\Im (J-B)^{-1} > 0$ and the map 
$B\mapsto   (J- B)^{-1}$ is a $2\times 2$ matrix-valued analytic function 
of $B \in \mathbb{H}_2^-$. By the uniqueness of analytic extension, 
we have $R_Y(B)=p (J- B)^{-1}$ for any $B\in\mathbb{H}_2^-$. 
An analogous argument holds for $B\in\mathbb{H}_2^+$. 
\end{proof}

\medskip

Using the above, and the fact that $\Omega$ is analytic it is easy to show that the mapping $B\mapsto B-R_Y(G_{X+Y}(B))$ satisfies the requirements for subordination as described in Subsection~\ref{subsec:2.2} hence we get the following lemma.
\begin{lemma}
 \label{lemma:22-Omega}
Suppose $y$ is a free Poisson element. 
For any $B\in \mathbb{H}_2^+$ we have  \begin{equation}\label{eqn:3.2.Omega1}\Omega(B)=B-R_Y(G_{X+Y}(B)).
\end{equation}
\end{lemma}

 For $B\in\mathbb{H}_2^+$ and $p>0$ we define
\begin{equation}
    H(B)=B+R_Y(G_X(B))=B+p(J-G_{X}(B))^{-1}.
\end{equation}
\hfill$\square$
\begin{lemma}
The identity \begin{equation}\label{eqn:H-left-inverse}
    H(\Omega(B))=B
\end{equation}
holds for arbitrary $B\in \mathbb{H}_2^+$.
\end{lemma}
\begin{proof}
Using $G_{X+Y}(B)=G_X(\Omega(B))$ we can rewrite \eqref{eqn:3.2.Omega1}  as
\[
  \Omega(B)=B-R_Y(G_X(\Omega(B)))
\]    
and the required result follows for $B\in\mathbb{H}_2^+$.
\end{proof}

\bigskip

\section{Properties of the functions $\delta_0$ and $h$}  \label{section:3}

\subsection{Proof of Proposition~\ref{prop:12}}
\label{sec:3.1}

\begin{lemma}  \label{lemma:integral-rep}
Let $\lambda\in\mathbb{C}$ and let $v\in\bC^2$ be a unit vector. 
There exist a constant $C$ and a positive Borel measure $\rho$ on $\mathbb{R}$ 
(depending on $\lambda$ and $v$), with $\rho$ satisfying the growth condition 
$\int_\mathbb{R}\frac{1}{1+t^2}d\rho(t)<\infty$,
such that the integral representation
\begin{equation}
   \label{eqn:integral-rep}
      \langle  \, (J-G_X(B(\lambda,\zeta)))^{-1}v \, , \,  v \, \rangle 
      = C + \int_\mathbb{R}\left(\frac{1}{\zeta-t}+\frac{t}{1+t^2}\right)d\rho(t),
      \ \ \zeta \in \bC^{+}
\end{equation}
holds.
\end{lemma}

\begin{proof}
Let $f : \bC^{+} \to \bC$ be defined by
\[
  f (\zeta):= 
  \langle \, (J-G_X(B(\lambda,\zeta)))^{-1}v \, , \,  v \, \rangle, \ \ \zeta \in \bC^{+}.
\]
The representation given for $f( \zeta )$ on the right-hand side of (\ref{eqn:integral-rep}) 
comes out as an application of the Nevanlinna integral representation theorem 
(see \cite[Section 5.3]{RRbook1994}), if we can prove that $f$ has the following 
properties: it maps $\bC^{+}$ to $\bC^{-}$, it is holomorphic, and satisfies 
the condition that 
$\lim_{\ee \rightarrow+ \infty} f( \qi\varepsilon) / \qi \ee = 0$.

Towards verifying the aforementioned properties of $f$, let us
denote by $G_{X,ij}(B)$ the $(ij)-$th entry of $G_X(B)$. Then we have 
\begin{align*}
 (J-G_X(B))^{-1}&=\begin{bmatrix}
   -G_{X,11}(B) &  \qi-G_{X,12}(B)\\
   -\qi-G_{X,21}(B) & -G_{X,22}(B)
 \end{bmatrix}^{-1}\\
 &=\frac{1}{\det(\lambda,\zeta)}\begin{bmatrix}
  -G_{X,22}(B) & -\qi+G_{X,12}(B)\\
  \qi+G_{X,21}(B) & -G_{X,11}(B)
 \end{bmatrix}
\end{align*}
where 
\[
 \det(\lambda,\zeta)=G_{X,11}(B)\cdot G_{X,22}(B)+(\qi-G_{X,12}(B))(\qi +G_{X,21}(B)).
\]
We note that for $B=B(\lambda,\zeta)$ with $\zeta\in\bC^+$ 
we have $\Im (J - G_X(B)) = - \Im G_X(B) >0$, which implies that
$\Im(J-G_X(B))^{-1}<0$. It follows that
\begin{align*}
  \Im f (\zeta)
  & = \Im  \langle \, (J-G_X(B(\lambda,\zeta)))^{-1}v \, , \,  v \, \rangle \\
  & =
   \langle \, \Im (J-G_X(B(\lambda,\zeta)))^{-1}v \, , \,  v \, \rangle < 0,
\end{align*}
hence that $f$ maps $\bC^{+}$ to $\bC^{-}$.  
The fact that $f$ is holomorphic on $\mathbb{C}^+$ holds because 
each of the entries $G_{X,ij}(B(\lambda,\zeta))$ is a holomorphic function on $\mathbb{C}^+$.

Finally, for the $\ee \to \infty$ limit, one first 
checks directly that $\lim_{\zeta\rightarrow\infty}G_{ij}(B(\lambda,\zeta))=0$ for $i,j\in\{1,2\}$,
which has the consequence that 
\[
\lim_{\zeta\rightarrow\infty}\det(\lambda,\zeta)=-1,
\]
and implies that
\[
  \lim_{\varepsilon\rightarrow+\infty}\frac{f(\qi\varepsilon)}{\qi\varepsilon}
  =\lim_{\varepsilon\rightarrow+\infty}
  \frac{\langle \, 
  ((J-G_X(B(\lambda,\qi \varepsilon)))^{-1}) v \, , \, v \, \rangle}{\qi\varepsilon}=0,
\]
as required.
\end{proof}

\vspace{6pt}

We are now ready to prove Proposition~\ref{prop:12}. 

\begin{proof}[Proof of Proposition~\ref{prop:12}]
$1^o$  We will use Lemma \ref{lemma:integral-rep} in the special case when the vector $v$
under consideration is 
$v=\begin{bmatrix}
    1 & 0
\end{bmatrix}^{\mathrm{T}} \in \bC^2$.
For this choice of $v$, the left-hand side of Equation (\ref{eqn:integral-rep}) is
readily seen to be $\bigl( H_{1,1} ( \lambda, \zeta ) - \zeta )/p$ (where recall that 
$H_{1,1} ( \lambda, \zeta ) := [ \, H ( \, B( \lambda , \zeta ) \, ) ]_{1,1}$, with 
$H(B) = B + \pp \,  \bigl( J - G_{X}(B) \bigr)^{-1}$ as defined in (\ref{eqn:H_Definition}). 
Moreover, for this choice of $v$, the positive measure $\rho$ provided by Lemma 
\ref{lemma:integral-rep} turns out to be a probability measure, and Equation 
(\ref{eqn:integral-rep}) takes the form
\begin{equation}    \label{eqn:31x}
H_{11}\left( B(\lambda,\zeta)  \right) - \zeta
= p \cdot \int_\mathbb{R}\frac{1}{\zeta-t}d\rho(t).
\end{equation}
For the justification of (\ref{eqn:31x}) we use the intermediate equality
\[
H_{11}\left( B(\lambda,\zeta)  \right) - \zeta
=-\frac{G_{X,22}(B)}{\text{det}(\lambda,\zeta)}
\]
(which was noticed during the proof of Lemma \ref{lemma:integral-rep}), and we invoke the fact that
$\lim_{\delta\rightarrow+\infty}\qi \delta \cdot G_{X,22}(B(\lambda,\qi\delta))=1$ and $\lim_{\delta\rightarrow\infty}\det(\lambda,\qi\delta)=-1$.

\vspace{6pt}

\noindent
$2^o$  We note that 
\[
  G_{X,11}(B(\lambda,\qi \delta))=G_{X,22}(B(\lambda,\qi \delta))
   =-\qi\delta \tau ((|\lambda-x|^2+\delta^2)^{-1}),
\]
and $G_{X,21}(B(\lambda,\qi \delta))=\overline{G_{X,12}(B(\lambda,\qi \delta))}$. It follows that 
$\Re H_{11}(B(\lambda,\qi \delta))=0$ for any $\delta\in (0,\infty)$. 

\vspace{6pt}

\noindent
$3^o$
For $\delta>0$ we have 
\[
   \frac{\Im H_{11}(B(\lambda,\qi \delta))}{\delta}
     =1-p\cdot\int_\mathbb{R}\frac{1}{\delta^2+t^2}d\rho(t)
\]
which is a stricly increasing function of $\delta\in (0,\infty)$. Moreover,
\[
 \lim_{\delta\rightarrow\infty} \frac{\Im H_{11}(B(\lambda,\qi \delta))}{\delta} =1. 
\]
\end{proof}

\subsection{Extension of $\delta_0$ and $h$.}\label{subsec:extension}
\phantom{x}

\noindent
Continuing from above, let us now recall that
Definition~\ref{def:13} introduced an open set $\cD \subseteq \bC$, where for every 
$\lambda\in\cD$ there exists a unique $\delta_0(\lambda)>0$  such that
\begin{equation}
  \label{eqn:identity-recall}
     \frac{\Im H_{11} \bigl( \, \lambda, \qi \delta_0(\lambda) \, \bigr)}{\delta_0(\lambda)} = 0.
\end{equation}

\vspace{6pt}

\begin{notation}
We let $\delta_0 : \bC \to [0, \infty)$ where for $\lambda \in \cD$ the value of $\delta_0 ( \lambda )$
is defined via the requirement in (\ref{eqn:identity-recall}), while for  
$\lambda\in \mathbb{C}\backslash\cD$ we put $\delta_0(\lambda) :=0$.
\end{notation}

\vspace{6pt}

\begin{proposition}
 \label{prop:delta-continuity}
The function $\lambda\mapsto\delta_0(\lambda)$ is  continuous  on $\mathbb{C}$. 
Moreover, it is differentiable in the open set $\cD$. 
\end{proposition}
\begin{proof}
The fact that the function $\lambda\mapsto\delta_0(\lambda)$ is differentiable on $\cD$ follows
from Proposition~\ref{prop:12} and the implicit function theorem. 

In order to establish the continuity on $\bC$ stated in the proposition, we are left to prove 
the following fact: if $\lambda$ is on the boundary of $\cD$ and if $\{\lambda_n\}$ is a sequence 
in $\cD$ such that $\lambda_n\rightarrow\lambda$, then it follows that 
    \begin{equation}\label{eqn:limit-delta}
\lim_{n\rightarrow\infty}\delta_0(\lambda_n) = 0.
    \end{equation}
    For the remaining part of the proof we fix such $\lambda$ and $\{\lambda_n\}$, for which we will 
    prove that (\ref{eqn:limit-delta}) holds.  
    
We proceed by contradiction.  Suppose that (\ref{eqn:limit-delta}) does not hold.
Dropping to a subsequence if necessary, we may assume that there exists $\varepsilon_0>0$ such that $\delta_0(\lambda_n)>\varepsilon_0$ for every $n$.
Then by the monotonicity result in Proposition~\ref{prop:12}, we have 
    \begin{equation}
       \label{eqn:inequality-ratio}
          \frac{\Im H_{11}(B(\lambda_n,\qi \delta_0(\lambda_n)))}{\delta_0(\lambda_n)}
      >\frac{\Im H_{11}(B(\lambda_n,\qi\varepsilon_0))}{\varepsilon_0}.
    \end{equation}
By the definition of $H$ in \eqref{eqn:H_Definition}, since $\varepsilon_0>0$, we deduce that 
\[
  \lim_{n\rightarrow\infty}H_{11}(B(\lambda_n,\qi \varepsilon_0))=H_{11}(B(\lambda,\qi \varepsilon_0)).
\]
Then \eqref{eqn:inequality-ratio} and Proposition \ref{prop:12} imply that 
\[
   \lim_{\varepsilon\rightarrow 0}\frac{\Im H_{11}(B(\lambda,\qi\varepsilon))}{\varepsilon}   
    <\frac{\Im H_{11}(B(\lambda,\qi\varepsilon_0))}{\varepsilon_0}\leq 0.
\]
This contradicts our assumption that $\lambda$ is on the boundary on $\cD$, and thus 
completes the proof.
\end{proof}

\begin{corollary}\label{cor:h-differentiable-on-D}
The function $h$ is differentiable on $\cD$. 
\end{corollary}
\begin{proof}
    The function $h(\lambda)=H_{1,2}(\lambda,\delta_0(\lambda))$ is a composition of two 
    differentiable functions, and as such it is differentiable.
\end{proof}

\vspace{6pt}

The next item to look at is the function $h : \cD \to \cM$, which was also introduced 
in Definition \ref{def:13}.  In the next remark we observe that under 
Assumption~\ref{assumption:16}
the formula (\ref{eqn:13c}) that defined $h$ can be extended to 
$\bC \setminus  (\cD\cup\{\alpha_1, \alpha_2, \cdots, \alpha_k\}) $. In the following we show that $h$ can be extended to $\bC$.
\vspace{6pt}

\begin{remark-and-definition}   \label{rem:34}
For $\lambda \notin \spec(x)$, the matrix 
$B( \lambda , 0 ) - X$ is invertible in $M_2 ( \cA )$ (where 
$B( \lambda , 0 ) := \left[   \begin{array}{cc}
0  & \lambda  \\ \overline{\lambda}  &  0 
\end{array}  \right]$ and 
$X := \left[   \begin{array}{cc}
0  & x  \\ x  &  0    \end{array}  \right]$).
The Cauchy transform $G_X$ is defined at $B( \lambda , 0)$, with
\begin{equation}    \label{eqn:34a}
G_X \bigl( \, B( \lambda , 0 ) \, \bigr)
= \mathbb{E} \left( \begin{bmatrix}
    0 & (\overline{\lambda}-x)^{-1}\\
    (\lambda-x)^{-1} & 0 
\end{bmatrix}\right)
      =\begin{bmatrix}
         0 & G_{\mu_x}(\overline{\lambda})\\
         G_{\mu_x}(\lambda) & 0
      \end{bmatrix}.
\end{equation}

Denote by $\{\alpha_1, \alpha_2, \cdots, \alpha_k\}$ the atoms of $\mu_x$ that are not in $\mathrm{cl}(\cD)$.
Under Assumption~\ref{assumption:16}, one immediately sees that
\eqref{eqn:34a} holds in $\bC \setminus  (\cD\cup\{\alpha_1, \alpha_2, \cdots, \alpha_k\})$.
From the explicit form of $G_X \bigl( \, B( \lambda , 0 ) \, \bigr)$ found in (\ref{eqn:34a}), 
one can check that $\qi +G_{\mu_x}(\lambda) \neq 0$ (see \eqref{eqn:inequality-outsidep}) 
and hence that $J - G_X \bigl( \, B( \lambda , 0 ) \, \bigr)$ is invertible
(where recall that $J := \begin{bmatrix} 0 & \qi\\ -\qi & 0 \end{bmatrix}$).  This allows
us to meaningfully plug $B = B( \lambda , 0)$ into the definition we had made 
for the function $H$ in (\ref{eqn:H_Definition}):
\[
H \bigl( \, B( \lambda , 0 ) \, \bigr) :=
B( \lambda, 0 ) + p \, 
\Bigl( J - G_X \bigl( \, B( \lambda , 0 ) \, \bigr) \, \Bigr)^{-1}.
 \]
Following to this, we can then define (in the same vein as we had done for 
elements of $\cD$ in Equation (\ref{eqn:13c}), but where now we have 
$\delta_0 ( \lambda ) = 0$):
\begin{equation}   \label{eqn:34c}
h( \lambda ) := \Bigl[ \, H \bigl( \, B( \lambda , 0 ) \, \bigr) \, \Bigr]_{1,2}, \qquad \lambda\in \bC \setminus  (\cD\cup\{\alpha_1, \alpha_2, \cdots, \alpha_k\}).
\end{equation}
\end{remark-and-definition}

\vspace{6pt}

\begin{notation-and-remark}   \label{def:41}
Consider some numbers $\alpha, \beta \in \bR$ and $\delta \in (0, \infty )$.  
It will be convenient to express the upcoming formulas in terms of
the following functions:
\begin{align}\label{eqn:Sdefinition}
    S&=S(\alpha,\beta,\delta) :=\int_{\bR} \frac{\alpha-t}{(\alpha-t)^2+\beta^2+\delta^2}\,d\mu_x(t),
\\
\label{eqn:Tdefinition}
    T&=T(\alpha,\beta,\delta) :=\int_{\bR} \frac{1}{(\alpha-t)^2+\beta^2+\delta^2}\,d\mu_x(t),
\\
\label{eqn:Ddefinition}
 D&=D(\alpha,\beta,\delta) :=\delta^2 T^2 + S^2 + (1-\beta T)^2.
\end{align}

An easy calculation shows that the Cauchy transform 
of $X$ at $B= B( \alpha + \qi \beta, \qi \delta ) = \begin{bmatrix}
    \qi \delta& \alpha+\qi \beta\\ \alpha-\qi \beta&\qi\delta
\end{bmatrix}$ is
$G_X (B) =
\begin{bmatrix}-\qi\delta T & S+\qi\beta T\\
        S - \qi\beta T & -\qi\delta T        \end{bmatrix},$
with the further consequence that 
\begin{align*}
    (J-G_X(B))^{-1}&=\begin{bmatrix}
  \qi\delta T&-S+\qi(1-\beta T)\\
  -S-\qi(1-\beta T)&\qi\delta T.
    \end{bmatrix}^{-1}
    \\&=-\frac{1}{D}\begin{bmatrix}
  \qi\delta T&S-\qi(1-\beta T)\\
  S+\qi(1-\beta T)&\qi\delta T.
    \end{bmatrix}.
\end{align*}
Recalling that $H(B)=B+p(J-G_X(B))^{-1}$, we thus get
\begin{align}    \label{eqn:H(B)}
    H(B)=\begin{bmatrix}
        \qi\delta\left(1-p\frac{T}{D}\right)
        &\alpha-p\frac{S}{D}+\qi\left(\beta+p\frac{1-\beta T}{D}\right)\\
        \alpha-p\frac{S}{D}-\qi\left(\beta+p\frac{1-\beta T}{D}\right)& \qi\delta\left(1-p\frac{T}{D}\right)
    \end{bmatrix} .
\end{align}
By comparing with Remark \ref{remark:1.9} and Definition \ref{def:13}, we note that  
\begin{equation}  \label{eqn:domain-D-new-def}
 \cD=\left\{\lambda=\alpha+\qi \beta\in\bC \middle\vert \lim_{\delta\to 0^+} 
 \frac{D(\alpha,\beta,\delta)}{ T(\alpha,\beta,\delta)}< p\right\}
\end{equation}
and that for $\lambda=\alpha+\qi \beta \notin \spec(x)$ we have
\begin{equation}
  \label{eqn:limit-T/D-formula}
      \lim_{\delta\rightarrow{ 0}} \frac{T(\alpha,\beta,\delta)}{D(\alpha,\beta,\delta)}
      = \frac{\varphi \bigl( \,  |\lambda-x|^{-2} \, \bigr) }{ |i+G_{\mu_x}(\lambda)|^2}
\end{equation}
\end{notation-and-remark} 

\medskip

\begin{lemma}
\label{lemma:subharmonic-1/F}
    The function $\lambda\mapsto \frac{\varphi (  |\lambda-x|^{-2} ) }{ |i+G_{\mu_x}(\lambda)|^2}$ is subharmonic in $\mathbb{C}\setminus\spec(x)$.
\end{lemma}
\begin{proof} This follows by writing
\[
\frac{ \varphi( |\lambda-x|^{-2} ) }{|i+G_{\mu_x}(\lambda)|^2}
 =\varphi\left( \frac{1}{|(i+G_{\mu_x}(\lambda))(\lambda-x)|^2} \right)
 =\int_\mathbb{R}\frac{1}{|(i+G_{\mu_x}(\lambda))(\lambda-t)|^2}d\mu_x(t),
\]
and by observing that $\lambda\rightarrow \frac{1}{(i+G_{\mu_x}(\lambda))(\lambda-t)}$ is 
complex analytic in $\mathbb{C}\setminus \text{spec}(x)$. 
\end{proof}

\medskip

\begin{proposition}    \label{prop:atoms}
Suppose that Assumption \ref{assumption:16} holds.  

\vspace{6pt}

\noindent
$1^o$ The boundary of $\cD$ can be characterized as  
\begin{align}  \label{eqn:boundary-D}
      \bd(\cD)&=\left\{\lambda=\alpha+\qi \beta\in\bC \middle\vert \lim_{\delta\to 0^+} \frac{D(\alpha,\beta,\delta)}{ T(\alpha,\beta,\delta)}= p\right\}\nonumber\\
      &=\left\{\lambda\in\bC \middle\vert \frac{|i+G_x(\lambda)|^2}{\varphi\left( |\lambda-x|^{-2} \right)} =p \right\}.
\end{align}
     
\vspace{6pt}

\noindent
$2^o$ Let $\alpha$ be an atom of $\mu_x$. Then we have 
   \[
    \lim_{\delta\rightarrow {0^+}}\frac{ D(\alpha,0,\delta)}{T(\alpha,0,\delta)}=\mu_x(\{\alpha\}).
   \]
Moreover, we have that
\[
(*) \hspace{1cm}
\left\{  \begin{array}{l}
\bigl[ \, \alpha \in \cD \, \bigl] \ \Rightarrow \ \bigl[ \, \mu_x (\{\alpha\})<p \, \bigr]  \\
                                                 \\
\bigl[ \, \alpha\in \bC\setminus\mathrm{cl}(\cD) \, \bigr] \ \Rightarrow 
\ \bigl[ \, \mu_x(\{\alpha\})>p \, \bigr] .
\end{array}  \right.
\]
[The third possibility that could be considered in $(*)$, that 
$\alpha \in \mathrm{cl} ( \cD ) \setminus \cD = \mathrm{bd} ( \cD )$, is ruled out by our
Assumption \ref{assumption:16}.]
\end{proposition}
     
\begin{proof}
If $\alpha+\qi\beta\in \mathbb{C}\backslash \spec(x)$, then the function 
$\frac{T(\alpha,\beta,\delta)}{D(\alpha,\beta,\delta)}$ is continuous in some neighborhood of 
$(\alpha, \beta, 0)$ and subharmonic by using \eqref{eqn:limit-T/D-formula} and Lemma \ref{lemma:subharmonic-1/F}. 
If $\alpha+\qi\beta\in \mathbb{C}\backslash (\cD\cup\spec(x))$, then 
$\lim_{\delta\to 0^+} \frac{D(\alpha,\beta,\delta)}{ T(\alpha,\beta,\delta)}\geq p$, 
with equality when $\alpha+\qi \beta\in \bd(\cD)$. 

Recall that Assumption \ref{assumption:16} asks for the support of the continuous part of 
$\mu_x$ to be contained in $\cD$. 
Because of this and of the considerations from the preceding paragraph, the boundary 
characterization stated in $1^o$ will follow once we prove the properties for atoms of $\mu_x$ 
that were stated in $2^o$.

We begin with a basic fact concerning the Cauchy transform $G_{\mu_x}$, that 
for every $t\in \bR$ one has
    \[
    \sphericalangle\lim_{z\to t}(z-t)G_{\mu_x}(z)=\mu_{x}(\{t\}).
    \]
Taking in the above limit $z=t+\qi\delta$ gives us that 
    $G_{\mu_x}(t+\qi \delta)=S(t,0,\delta)-i\delta T(t,0,\delta)$, hence
    \[
    \lim_{\delta\to 0^+} \delta \qi (S(t,0,\delta)-\delta  T(t,0,\delta)\qi)=\mu(\{t\}),
    \]
and hence
    \begin{align}
      \label{limt-S-T}
        \lim_{\delta\to 0^+} \delta S(t,0,\delta)=0\quad
        \text{and}\quad 
        \lim_{\delta\to 0^+} \delta^2 T(t,0,\delta)=\mu_x(\{t\}).
    \end{align}
From the second limit in (\ref{limt-S-T})
it follows in particular that if $t\in\bR$ is an atom for $\mu_x$, 
then $\lim_{\delta\to 0^+} T(t,0,\delta)=+\infty$. 

Our next observation is that, for any atom $\alpha$ of $\mu_x$, one has that
    \begin{equation}
      \label{eqn:D-over-T-general}
        \lim_{\delta\rightarrow {0^+}}\frac{ D(\alpha,0,\delta)}{T(\alpha,0,\delta)}
        = \lim_{\delta\to 0^+}\left(\frac{1}{T(\alpha,0,\delta)}+\frac{S(\alpha,0,\delta)^2}{T(\alpha,0,\delta)}
        +\delta^2 T(\alpha,0,\delta)\right)=\mu_x(\{\alpha\}).
    \end{equation} 
Indeed, the first term $1/ T(\alpha,0,\delta)$
listed in the above limit converges to $0$ for $\delta \to 0^{+}$, because  
$T(\alpha,0,\delta)$ converges to $\infty$. Then for the second term listed in
(\ref{eqn:D-over-T-general}), we have
     \[
       \lim_{\delta\to 0^+}\frac{S(\alpha,0,\delta)^2}{T(\alpha,0,\delta)}=
        \lim_{\delta\to 0^+}\frac{ ( \delta S(\alpha,0,\delta) )^2}{\delta^2 T(\alpha,0,\delta)}=0,
     \]
where once again we invoke \eqref{limt-S-T}.  Finally, the third term
$\delta^2 T(\alpha,0,\delta )$ is precisely the one giving the mass $\mu_x (\{ \alpha \})$.

As a further consequence of (\ref{eqn:D-over-T-general}), 
combined with (\ref{eqn:domain-D-new-def}),
we find that for an atom $\alpha$ of $\mu_x$ we have
\begin{equation}   \label{eqn:38x}
\alpha \in \cD  \ \Leftrightarrow \
\lim_{\delta\rightarrow {0^+}}\frac{ D(\alpha,0,\delta)}{T(\alpha,0,\delta)} < p
\ \Leftrightarrow \ \mu_x ( \{ \alpha \} ) < p.
\end{equation}
In particular, this gives the first implication $(*)$ in the statement of part $2^o$.

By Assumption \ref{assumption:16}, an element in the interior set 
$\mathrm{int}(\mathbb{C}\setminus \cD)=\mathbb{C}\setminus \mathrm{cl}(\cD)$ either is not in $\spec(x)$ at all, 
or is an atom of $\mu_x$. If $\alpha_0\in \mathbb{C}\setminus \mathrm{cl}(\cD)$ is an atom of $\mu_x$, then we 
first note that, due to the equivalence (\ref{eqn:38x}) obtained above, we have $\mu_x(\{\alpha_0\})\geq {p}$. Moreover 
$\mu_x(\{\alpha_0\})\neq p$ due to the fact that $\lim_{\delta\rightarrow{0}}\frac{T(\alpha,\beta,\delta)}{D(\alpha,\beta,\delta)}\leq \frac{1}{p}$ for $\alpha+\qi\beta\in\mathbb{C}\setminus \mathrm{cl}(\cD)$, the identity \eqref{eqn:limit-T/D-formula}, and 
Lemma \ref{lemma:subharmonic-1/F}. Hence $\mu_x(\{\alpha_0\})>p$, which establishes 
the second implication $(*)$ in the statement of part $2^o$.
\end{proof}

\medskip

\begin{proposition}  
\label{prop:outside-domain-injective}
\noindent
Suppose that Assumption~\ref{assumption:16} holds.  Then the set of atoms of $\mu_x$
that are not in $\mathrm{cl}(\cD)$ is finite, say $\{\alpha_1, \alpha_2,\cdots, \alpha_k\}$, 
and we can consider the function 
$h : \bC \setminus  (\cD\cup\{\alpha_1, \alpha_2, \cdots, \alpha_k\})\to \bC$
defined in \eqref{eqn:34c}. 

\vspace{6pt}

\noindent
$1^o$ $h$ is analytic on $\bC\setminus (\mathrm{cl}(\cD)\cup\{\alpha_1, \alpha_2, \cdots, \alpha_k\})$, 
where one has 
\[
h(\lambda)=\lambda-\frac{p}{\qi+G_{\mu_x}(\lambda)}, 
\mbox{ with } G_{\mu_x}(\lambda) := \int_\mathbb{R}\frac{1}{\lambda-u}d\mu_x(u).
\]
Moreover, $h$ has removable singularities at $\alpha_1, \alpha_2, \cdots, \alpha_k$,
which allows us to extend it analytically on
$\bC \setminus \mathrm{cl} ( \cD )$, by putting $h(\alpha_i)=\alpha_i$ for $1 \leq i \leq k$.
In addition, we have  $0<h'(\alpha_i)<1$ and $\mu_x(\{\alpha_i\})>p$.

\vspace{6pt}

\noindent
$2^o$ One has that
\[
0 < |h(\lambda_1)-h(\lambda_2)| \leq 2 |\lambda_1-\lambda_2|, \ \forall \, \lambda_1 \neq \lambda_2
\ \mathrm{ in } \ \bC \setminus \cD.
\]
In particular, $h$ is injective on $\mathbb{C}\setminus \cD$. 

\end{proposition}

\begin{proof}
For $\lambda\in \bC\setminus (\cD\cup\{\alpha_1, \alpha_2, \cdots, \alpha_k\})$ we have
that $\delta_0(\lambda)=0$, and (\ref{eqn:34a}) follows from Assumption~\ref{assumption:16},
since $\spec(x)\subset \cD\cup\{\alpha_1, \alpha_2, \cdots, \alpha_k\}$.
    Recall that $H(b)=b+p(J-G_X(b))^{-1}$ where $J=\begin{bmatrix}
      0 & \qi\\
      -\qi & 0
    \end{bmatrix}$. 
 Hence, 
    \begin{align*}
     h(\lambda)&=H_{12}\left(\begin{bmatrix}
         0 & \lambda\\
         \overline{\lambda} & 0
\end{bmatrix}             \right)
    =\lambda-\frac{p}{\qi+G_{\mu_x}(\lambda)},
    \end{align*}
which is analytic on $\bC\setminus (\mathrm{cl}(\cD)\cup\{\alpha_1, \alpha_2, \cdots, \alpha_k\})$. Moreover,
\[
   \lim_{\lambda\rightarrow{\alpha_i}}h(\lambda)=\alpha_i.
\]
Hence, $h$ has removable singularities at $\{\alpha_1, \alpha_2, \cdots, \alpha_k\}$ such that \[h(\alpha_i)=\alpha_i\] for $i=1,\cdots, k$. 
Moreover, we can calculate the derivative as follows:
\begin{align*}
    h'(\alpha_i)=1-\lim_{\lambda\rightarrow{\alpha_i}}\frac{p}{(\qi +G_{\mu_x}(\alpha_i))(\lambda-\alpha_i)
    }
    =1-\frac{p}{\mu_x(\alpha_i)}.
\end{align*}
By Proposition \ref{prop:atoms}, we deduce that $\mu_x(\{\alpha_i\})>p$ and hence $0<h'(\alpha_i)<1$.

For $\lambda\in \bC\setminus (\cD\cup\{\alpha_1, \alpha_2, \cdots, \alpha_k\})$, by Remark \ref{remark:1.9}, we have 
\begin{equation}
\label{eqn:3.17}
    \varphi\left( \,|\lambda-x|^{-2} \, \right)
  \cdot \frac{1}{|\qi+ G_{\mu_x}(\lambda)|^2}=
  \int_\mathbb{R}\frac{1}{|\lambda-u|^2}d\mu_x(u)\frac{1}{|\qi+ G_{\mu_x}(\lambda)|^2}\leq \frac{1}{p}. 
\end{equation}
Suppose $\lambda_1, \lambda_2\in \bC\setminus (\cD\cup\{\alpha_1, \alpha_2, \cdots, \alpha_k\})$ and $\lambda_1	\neq \lambda_2$. We have 
\begin{align}
  h(\lambda_1)-h(\lambda_2)&=(\lambda_1-\lambda_2)+p\frac{G_{\mu_x}(\lambda_1)-G_{\mu_x}(\lambda_2)}{(\qi +G_{\mu_x}(\lambda_1))(\qi +G_{\mu_x}(\lambda_2))}\nonumber\\
   &=(\lambda_1-\lambda_2)p\left(\frac{1}{p}-\frac{\int_\mathbb{R}\frac{1}{(\lambda_1-u)(\lambda_2-u)}d\mu_x(u)}{{(\qi +G_{\mu_x}(\lambda_1))(\qi +G_{\mu_x}(\lambda_2))}} \right).\label{eqn:3.18}
\end{align}
By Cauchy-Schwarz inequality, we have 
\begin{align*}
  \left|\int_\mathbb{R}\frac{1}{(\lambda_1-u)(\lambda_2-u)}d\mu_x(u) \right|^2
  \leq \left(\int_\mathbb{R}\frac{1}{|\lambda_1-u|^2}d\mu_x(u)  \right)\left(\int_\mathbb{R}\frac{1}{|\lambda_2-u|^2}d\mu_x(u)  \right).
\end{align*}
Denote 
\[
K(\lambda_1,\lambda_2)= \left(\int_\mathbb{R}\frac{1}{|\lambda_1-u|^2}d\mu_x(u)  \right)\left(\int_\mathbb{R}\frac{1}{|\lambda_2-u|^2}d\mu_x(u)  \right)-\left|\int_\mathbb{R}\frac{1}{(\lambda_1-u)(\lambda_2-u)}d\mu_x(u) \right|^2.
\]
 For arbitrary $f,g\in L^2(\mu_x)$, we have 
\begin{equation}\label{eq:lagrange}
\|f\|_2^2\|g\|_2^2-|\langle f,\overline{g}\rangle|^2
=
\frac12\iint |f(u)g(v)-f(v)g(u)|^2 d\mu_x(u) d\mu_x(v).
\end{equation}
We then apply this identity for $f(u)=\frac{1}{\lambda_1-u}$ and 
$g(u)=\frac{1}{\lambda_2-u}$. 
One can check 
\[
\frac{1}{(\lambda_1-u)(\lambda_2-v)}
  -\frac{1}{(\lambda_1-v)(\lambda_2-u)}
=\frac{(\lambda_1-\lambda_2)(v-u)}
       {(\lambda_1-u)(\lambda_1-v)
        (\lambda_2-u)(\lambda_2-v)}.
\]
Hence we have 
\begin{equation}\label{eq:defect-identity}
\begin{aligned}
K(\lambda_1,\lambda_2)
  =\frac{|\lambda_1-\lambda_2|^2}{2}
    \iint_{\mathbb{R}^2}
    \frac{|u-v|^2}
    {|\lambda_1-u|^2|\lambda_1-v|^2
     |\lambda_2-u|^2|\lambda_2-v|^2}
    d\mu_x(u)d\mu_x(v).
\end{aligned}
\end{equation}
Therefore, for $\lambda_1\neq\lambda_2$ we have $K(\lambda_1,\lambda_2)\geq 0$ and the inequality 
is strict unless $\mu_x$ is concentrated in a single point. 

By combining with inequality \eqref{eqn:inequality-outsidep}, we obtain
\[
  \left| \frac{\int_\mathbb{R}\frac{1}{(\lambda_1-u)(\lambda_2-u)}d\mu_x(u)}{{(-\qi -G_{\mu_x}(\lambda_1))(-\qi -G_{\mu_x}(\lambda_2))}} \right| <\frac{1}{p}. 
\]
Therefore, by \eqref{eqn:3.18} we have $h(\lambda_1)\neq h(\lambda_2)$. 
The above estimate also shows that $|h(\lambda_1)-h(\lambda_2)|<2|\lambda_1-\lambda_2|$. 

It remains to show that if $\lambda\in \bC\setminus (\cD\cup\{\alpha_1, \alpha_2, \cdots, \alpha_k\})$ and $1\leq j\leq k$, we have $h(\lambda)\neq \alpha_i=h(\alpha_i)$. Suppose that $h(\lambda)=\alpha_i$, then we have 
\[
  (\lambda-\alpha_i)(\qi +G_{\mu_x}(\lambda))=p,
\]
hence $\vert (\lambda-\alpha_i)(\qi +G_{\mu_x}(\lambda)) \vert=p$.
Then \eqref{eqn:3.17} implies that 
\begin{align*}
   1\geq \varphi\left( \,|\lambda-x|^{-2} \, \right)
  \cdot \frac{p}{|\qi+ G_{\mu_x}(\lambda)|^2}>\frac{\mu_x(\{\alpha_i\})}{|\lambda-\alpha_i|^2}\frac{p}{|\qi +G_{\mu_x(\lambda)}|^2}=\frac{\mu_x(\{\alpha_i\})}{p},
\end{align*}
which contradicts the fact that $\mu_x(\{\alpha_i\})>p$ from Proposition \ref{prop:atoms}.
This completes the proof for $2^o$.
\end{proof}

\subsection{$H$ is the bilateral inverse of $\Omega$}
\phantom{x}

\noindent
The following result is an analogue of \cite[Lemma 4.2]{BelinschiCapitaine}. The proof is adapted from the argument there. 

\begin{proposition}\label{prop:3.2}
Fix $\lambda \in \bC$. If $\zeta \in \bC^{+}$ has the property that 
$\Im H(B(\lambda,\zeta))>0$, then it follows that 
$\Im H(B(\lambda,\zeta+\qi \varepsilon))>0$
for every $\varepsilon>0$. Moreover, for $B = B( \lambda, \zeta )$ with $\zeta$ 
as above, one has that
\begin{equation}\label{eqn:3.5.inverseEqn}
  \Omega_1(H(B))=B.
\end{equation}
and 
\begin{equation}
  G_{X+Y}(H(B))=G_X(B).
\end{equation}
\end{proposition}

\begin{proof}
Let $\zeta = \alpha_o + \qi \beta_o \in \bC^{+}$ be such that 
$\Im H(B(\lambda,\zeta))>0$.  Let us also consider a unit vector
$v \in \mathbb{C}^2$. By \eqref{eqn:integral-rep} in 
Lemma~\ref{lemma:integral-rep}, there exist a constant $C$ and a Borel measure 
$\rho$ such that
\[
 \langle  H(B(\lambda, \alpha_o + \qi \beta ))v,  v\rangle 
   =\left\langle B(\lambda, \alpha_o + \qi \beta )v ,v \right\rangle
    +C+
   p\cdot\int_\mathbb{R}\left(\frac{1}{\alpha_o + \qi \beta -t}
   +\frac{t}{1+t^2}\right)d\rho(t),
\]
for every $\beta > 0$.
Note that 
$\Im (\langle B(\lambda,\alpha_o + \qi \beta ) v, v \rangle) = \beta$. Hence, 
\[
  \Im \langle H(B(\lambda,\alpha_o + \qi \beta ))v, v \rangle 
  = \beta \left( 1 -p \cdot \int_\mathbb{R}\frac{1}{\beta^2+(\alpha_o - t)^2}d\rho(t) \right).
\]
It follows that the map 
$\beta \mapsto \frac{\Im (\langle  H(B(\lambda, \alpha_o +\qi \beta))v, v \rangle ) }{\beta}$ 
is increasing on $(0,\infty)$. 
Since this is true for any vector $v\in\mathbb{C}^2$, it follows that the positivity property 
given to hold at $\beta = \beta_o$ must also hold for $\beta > \beta_o$. This gives the 
positivity of $\Im H(B(\lambda,\zeta+\qi \varepsilon))$ for $\ee > 0$, since
$\zeta + \qi \ee$ is of the form $\alpha_o + \qi \beta$ with $\beta > \beta_o$.
  
We next remark that, for the same $\zeta$ as above and for $\varepsilon$ sufficiently large, 
the inverse function theorem yields the equality $\Omega[H(B(\lambda,\zeta+\qi \varepsilon))]=B(\lambda,\zeta+\qi \varepsilon)$. We note that the set $\{ \zeta\in\mathbb{C}^+: \Im H(B(\lambda,\zeta))>0\}$ is open. Moreover, each entry of $H(B(\lambda,\delta))$ is a holomorphic function of $\zeta$ in the entire upper half plane $\mathbb{C}^+$. By the identity principle for analytic functions, we conclude that
\[
   \Omega[H(B(\lambda,\zeta+\qi \varepsilon))]=B(\lambda,\zeta+\qi \varepsilon)
\]
for any $\varepsilon\geq 0$. In particular, $\Omega(H(B))=B$ for $B=B(\lambda,\zeta)$. Thus $\Im H(B(\lambda,\zeta))>0$ implies
\[
  G_{X+Y}(H(B))=G_{X}(\Omega(H(B)))=G_X(B).
\]
\end{proof}

\begin{corollary}\label{cor:imageSubord}
The image of the subordination function $\Omega$  can be described as follows
 \[
   \Omega\left( \left\{ B(z,\qi \varepsilon): z\in\mathbb{C}, \varepsilon>0 \right\} \right)   =
    \left\{ B(\lambda,\qi \delta): \lambda\in\mathbb{C}, \delta>\delta_0(\lambda) \right\}.
\]
\end{corollary}

\vspace{0.5cm}

\section{Discussion of the technical assumptions}  \label{section:4}

In this section we discuss Assumption~\ref{assumption:16}, and we show that it is satisfied 
in a number of examples.

\medskip

\begin{remark}   \label{rem:42}
We record here a further observation related to the special case $\beta = 0$ in Notation 
and Remark \ref{def:41}.  We have 
\[
S(\alpha,0,\delta)=\int_\mathbb{R}\frac{\alpha-t}{(\alpha-t)^2+\delta^2}d\mu_x(t), \quad
T(\alpha,0,\delta)=\int_\mathbb{R}\frac{1}{(\alpha-t)^2+\delta^2}d\mu_x(t).
\]
Observe that by the Cauchy-Schwarz inequality
\begin{equation}
\label{eq:CS:S2ltT}
    S(\alpha,0,\delta)^2\leq \int_\bR \frac{(\alpha-t)^2}{(\alpha-t)^2+\delta^2}d\mu_x(t)
    \cdot \int_\bR \frac{1}{(\alpha-t)^2+\delta^2}d\mu_x(t) \leq 1 \cdot T(\alpha, 0, \delta).
\end{equation}
We thus record the inequality
\begin{align}\label{eqn:ST_ineq}
\frac{S(\alpha,0,\delta)^2}{T(\alpha,0,\delta)}\leq 1,
\end{align}
where we also note that (\ref{eqn:ST_ineq}) holds with equality if and only if 
$\mu_x$ is concentrated at a point. 
\end{remark}

\medskip

We now take on the topic of the section and we first give sufficient conditions 
for Assumption~\ref{assumption:16} to hold.

\medskip

\begin{proposition}  \label{prop:43}
Assume that the distribution $\mu_x$ has the following properties:
    \begin{itemize}
        \item[$(i)$] There are no atoms of $\mu_x$ such that $\mu_x(\{\alpha\})=p$. 

        \item[$(ii)$] The limit  $G_{\mu_{x}}(\alpha):=\lim_{\beta\to 0^+}G_{\mu_{x}}(\alpha+i \beta)$ exists and is finite for every $\alpha\in\spec(x)$ which is not an atom of $\mu_{x}$.
        \item[$(iii)$] If for some $\alpha\in\spec(x)$ we have that $\Im G_{\mu_{x}}(\alpha)=0$, then $\lim_{\beta\to 0+}-\frac{\Im G_{\mu_{x}}(\alpha+i \beta)}{\beta^\gamma}= c>0$ for some $\gamma\in(0,1)$.
        \end{itemize}
        Then: \begin{itemize}
            \item[(1)] Assumption \ref{assumption:16} is satisfied. Moreover, by Proposition \ref{prop:atoms},
             $\spec(x)\subset \cD\cup\{\alpha_1,\ldots,\alpha_k\}$, where $\mu_x(\{\alpha_i\})>p$ for all $i=1, \cdots, k$.
            \item[(2)] The function 
        \begin{equation}
         \label{eqn:F-lambda-domain-D}
         F(\lambda) := \frac{ |i+G_x(\lambda)|^2 }{ \varphi \left( |\lambda-x|^{-2} \right) },
\quad \lambda\not\in\spec(x),
        \end{equation}
        can be extended continuously to $\spec(x)$ and (with $F$ denoting the extended function) one can equivalently describe $\cD$ as
        \[\cD=\{\lambda\in \bC:F(\lambda)<p\}.\]
        \end{itemize} 
\end{proposition}

\begin{remark}
    Point $(i)$ above shows that Assumption \ref{assumption:16} is satisfied for any purely discrete measure with 
    finitely many atoms and without atoms of weight equal to $p$. We point out that many measures 
    arising in free probability (e.g. semicircle distribution) are such that their Cauchy transform 
    has continuous continuation to the real line, so they satisfy assumption $(ii)$, moreover 
    \cite{BelinschiLebesgueDecomposition} shows that if $\sigma=\mu\boxplus\nu$, then $G_\sigma(z)$ 
    has continuous continuation to $\bR$, although the value $\infty$ is allowed, which we exclude above. 
    Assumption $(iii)$, roughly speaking, controls the behaviour of the density of $\mu_x$ at the edges
    of the support.
\end{remark}

\vspace{6pt}

\noindent
{\em Proof of Proposition \ref{prop:43}} We proceed similarly to the proof of Proposition \ref{prop:atoms}.

  \noindent
    (1) Pick $t \in \spec (x)$ which is not an atom for $\mu_x$, and prove that $t \in \cD$.
This amounts to proving that
     $\lim_{\delta\to 0^+}\frac{D(t,0, \delta)}{T(t,0, \delta)}< p$.

From assumption $(ii)$ made in the statement of the proposition it follows that
    $\lim_{\delta\to 0^+}S(t,0,\delta)$ and $\lim_{\delta\to 0^+}\delta T(t,0,\delta)$ exist and are finite.
 When $\lim_{\delta\to 0^+}\delta T(t,0,\delta)=0$ we still have $\lim_{\delta\to 0^+} T(t,0,\delta)=+\infty$ 
   and hence all three limits
\[
\lim_{\delta\to 0^+} \frac{1}{T(t,0,\delta)}, 
\ \lim_{\delta\to 0^+} \frac{S(t,0,\delta)^2}{T(t,0,\delta)} \mbox{ and }
\lim_{\delta\to 0^+} \delta^2 T(t,0,\delta)
\]
are equal to $0$. Hence $\lim_{\delta\rightarrow {0^+}}\frac{ D(\alpha,\beta,\delta)}{T(\alpha,\beta,\delta)}=0<p$, 
as we wanted.

\vspace{6pt}

    \noindent
    (2) Let $\lambda = \alpha + \qi\beta \in \bC \setminus \spec(x)$ with $\beta > 0$. We have
\begin{align}
\varphi\!\left(|\lambda-x|^{-2}\right)
 &= \int_{\bR} \frac{1}{|\alpha+\qi\beta-t|^2}\, d\mu_x(t) \notag\\
 &= \int_{\bR} \frac{1}{(\alpha-t)^2+\beta^2}\, d\mu_x(t)
 = T(\alpha,\beta,0),
 \label{eqn:phi_T}
 \end{align}
and
\begin{align}
G_{\mu_x}(\lambda)
 &= \int_{\bR} \frac{1}{\alpha+\qi\beta-t}\, d\mu_x(t) \notag\\
 &= \int_{\bR}
 \frac{\alpha-t-\qi\beta}{(\alpha-t)^2+\beta^2}\, d\mu_x(t)
 \notag\\
 &= S(\alpha,\beta,0)-\qi\beta T(\alpha,\beta,0).
 \label{eqn:Gx_ST}
\end{align}
Using \eqref{eqn:Gx_ST}, we can rewrite the numerator of the function $F(\lambda)$:
\begin{align*}
|\qi+G_{\mu_x}(\lambda)|^2
 &= \bigl|\qi + S(\alpha,\beta,0)-i\beta T(\alpha,\beta,0)\bigr|^2 \\
 &= \bigl|S(\alpha,\beta,0)
   + \qi\bigl(1-\beta T(\alpha,\beta,0)\bigr)\bigr|^2 \\
 &= S(\alpha,\beta,0)^2
   + \bigl(1-\beta T(\alpha,\beta,0)\bigr)^2.
\end{align*}
Substituting this and \eqref{eqn:phi_T} into the definition of $F(\lambda)$, we see that for $\lambda \notin \spec(x)$ the function evaluates exactly to the ratio $\frac{D}{T}$ at $\delta=0$:
\begin{equation}
\label{eqn:F_lambda_DT}
F(\alpha+\qi\beta)
=
\frac{
S(\alpha,\beta,0)^2
+
\bigl(1-\beta T(\alpha,\beta,0)\bigr)^2
}{
T(\alpha,\beta,0)
}
=
\frac{D(\alpha,\beta,0)}{T(\alpha,\beta,0)}.
\end{equation}
To extend $F$ continuously to the real axis at a point $\alpha \in \spec(x)$, we evaluate
\[
\lim_{\beta\to0^+} F(\alpha+\qi\beta).
\]
Observe that (\ref{eqn:F_lambda_DT}) can be rewritten as
\[
F(\alpha+\qi\beta) = \frac{1}{T(\alpha,\beta,0)}
+ \frac{S(\alpha,\beta,0)^2}{T(\alpha,\beta,0)}
+\beta^2T(\alpha,\beta,0)-2\beta,
\]
which is further processed as
\begin{align*}
F(\alpha+\qi\beta) + 2 \beta 
& = \frac{1}{T(\alpha,\beta,0)} + \frac{S(\alpha,\beta,0)^2}{T(\alpha,\beta,0)}
    + \beta^2T(\alpha,\beta,0)   \\
& = \frac{1}{T(\alpha,0, \beta)} + \frac{S(\alpha,0,\beta)^2}{T(\alpha,0,\beta)}
    + \beta^2T(\alpha,0,\beta)   \\
& = \frac{ D(\alpha,0,\beta) }{ T(\alpha,0,\beta) },
\end{align*}
where at the second equality sign we used the fact that $\beta$ and $\delta$ play symmetric roles in 
the definitions of integrals $S(\alpha,\beta,\delta)$ and $T(\alpha,\beta,\delta)$.

It then follows that $\lim_{\beta\to 0^+} F(\alpha+i\beta)$ exists. 
Moreover the continuous extension of $F(\lambda)$ to any
$\alpha\in\spec(x)$ (which we still denote by $F)$ has the same value as $\lim_{\delta\to 0^+}\frac{D(\alpha,0,\delta)}{T(\alpha,0,\delta)}$. Therefore, we can write
\[
\cD
=
\{\lambda\in\bC : F(\lambda)<p\},
\]
which completes the proof.
\hfill $\square$

\medskip

\begin{remark}
For $p <1$, one has $\mu_y(\{0\})=1-p$. So, for such $p$, by analogy to the selfadjoint case, 
one might expect that if $\mu_x$ has an atom $\alpha\in\bR$ with 
$\mu_x \bigl( \, \{ \alpha \} \, \bigr) + (1-p) > 1$ (that is, with
$\mu_x \bigl( \, \{ \alpha \} \, \bigr) > p$), then the Brown measure
$\mu_{x+\qi y}$ has an atom at $\alpha$, of mass $\mu_x \bigl( \, \{ \alpha \} \, \bigr) - p$.
On the other hand, if $p>\max_{t}\{\mu_x(\{t\}):t\in\mathbb{R}\}$ then we expect that 
$\mu_{x+\qi y}$ has no atom. 
These questions are discussed further in Section 7 below.
\end{remark}

\medskip

\section{Non-vanishing Jacobians}
\label{section:5}

In this section, we consider the mapping \(H\) on matrices of the form
$B(\alpha+\qi\beta, \qi\delta)=\begin{bmatrix}
\qi \delta & \alpha + \qi \beta\\[2mm]
\alpha - \qi \beta & \qi \delta
\end{bmatrix}$
where \(\alpha, \beta \in \mathbb{R}\) and \(\delta > 0\). As discussed before, the image of \(H\) has the same structure: the \((1,1)\) and \((2,2)\) entries coincide and are purely imaginary with positive imaginary part, while the \((1,2)\) and \((2,1)\) entries are complex conjugates of each other. This structural preservation allows us to view \(H\) as a mapping from \(\mathbb{R}^2 \otimes \mathbb{R}^+\) to itself, we denote this mapping by $\widehat{H}$.

Recall the definitions of $S$ and $T$ in \eqref{eqn:Sdefinition} and \eqref{eqn:Tdefinition}.
It is useful to record the derivatives of $S$ and $T$ with respect to $\alpha,\beta,\delta$.
\begin{proposition}\label{Prop:STderivatives} We have
    \begin{align*}
        \frac{\partial S}{\partial \alpha}&=T-2 \int_{\bR} \frac{(\alpha-t)^2}{\left((\alpha-t)^2+\beta^2+\delta^2\right)^2}d\mu_x(t),\\
        \frac{\partial S}{\partial \beta}&=-2\beta \int_{\bR} \frac{\alpha-t}{\left((a-t)^2+\beta^2+\delta^2\right)^2}d\mu_x(t)\\ 
        \frac{\partial S}{\partial \delta}&=-2\delta \int_{\bR} \frac{\alpha-t}{\left((\alpha-t)^2+\beta^2+\delta^2\right)^2}d\mu_x(t),\\
        \frac{\partial T}{\partial \alpha}&=-2 \int_{\bR} \frac{\alpha-t}{\left((\alpha-t)^2+\beta^2+\delta^2\right)^2}d\mu_x(t),\\
        \frac{\partial T}{\partial \beta}&=-2\beta \int_{\bR} \frac{1}{\left((\alpha-t)^2+\beta^2+\delta^2\right)^2}d\mu_x(t)\\ 
        \frac{\partial T}{\partial \delta}&=-2\delta \int_{\bR} \frac{1}{\left((\alpha-t)^2+\beta^2+\delta^2\right)^2}d\mu_x(t).
    \end{align*}
    Note that $\frac{\partial S}{\partial \beta}=\frac{\beta}{\delta}\frac{\partial S}{\partial \delta}$ and likewise $\frac{\partial T}{\partial \beta}=\frac{\beta}{\delta}\frac{\partial T}{\partial \delta}$.
\end{proposition}

\subsection{The Jacobian inside the domain $\cD$.}
\phantom{x}

\noindent
Recall that for $B=B(\alpha+\qi\beta,\qi\delta)$, we have 
\begin{align*}
    H(B)=\begin{bmatrix}
        \qi\delta\left(1-p\frac{T}{D}\right)
        &\alpha-p\frac{S}{D}+\qi\left(\beta+p\frac{1-\beta T}{D}\right)\\
        \alpha-p\frac{S}{D}-\qi\left(\beta+p\frac{1-\beta T}{D}\right)& \qi\delta\left(1-p\frac{T}{D}\right)
    \end{bmatrix},
\end{align*}
where $D=D(\alpha,\beta,\delta)=\delta^2T^2+S^2+(1-\beta T)^2$.

\vspace{6pt}

\begin{proposition}  \label{prop:JacobianH}
    For any $\alpha,\beta,\delta\in\bR$ denote 
    \begin{align*}
    \widehat{H}_{12}(\alpha,\beta,\delta)&=H_{12}\left(\begin{bmatrix}
        \qi\delta&\alpha+\qi\beta\\
        \alpha-\qi\beta& \qi\delta
    \end{bmatrix}\right) 
    \\
\widehat{H}_{11}(\alpha,\beta,\delta)
&= 
H_{11}\!\left(
\begin{bmatrix}
\qi\delta & \alpha + \qi\beta \\
\alpha - \qi\beta & \qi\delta
\end{bmatrix}
\right).
\end{align*}
Then for every $\alpha+\qi\beta\in\cD$ such that $\delta_0=\delta_0(\alpha+\qi\beta)>0$ the Jacobian of the mapping 
\begin{equation}   \label{eqn:JacobianH-1}
\widehat{H}:\begin{pmatrix}\alpha \\ \beta \\ \delta\end{pmatrix}
\to \begin{pmatrix}
    \Re\left(\widehat{H}_{12}\left(\alpha,\beta,\delta\right)\right)\\
    \Im\left(\widehat{H}_{12}\left(\alpha,\beta,\delta\right)\right)\\
    \Im\left(\widehat{H}_{11}\left(\alpha,\beta,\delta\right)\right)
\end{pmatrix}.
\end{equation} 
satisfies the inequality
    \[
    \det\left(J_{\widehat{H}}(\alpha,\beta,\delta_0)\right)>0.
    \]
\end{proposition}

\begin{proof}
Observe that the map $\widehat{H}$ defined in 
(\ref{eqn:JacobianH-1}) can be written as
\begin{equation}   \label{eqn:JacobianH-2}
    \widehat{H}(\alpha,\beta,\delta) =
    \begin{pmatrix}
    \alpha-p\frac{S}{(1-\beta T)^2+S^2+(\delta T)^2}\\ 
    \beta+p\frac{1-\beta T}{(1-\beta T)^2+S^2+(\delta T)^2}\\ 
    \delta-p\frac{\delta T}{(1-\beta T)^2+S^2+(\delta T)^2}
    \end{pmatrix}=
    \begin{pmatrix}
         \alpha-\frac{pS}{D}\\
         \beta+\frac{p(1-\beta T)}{D}\\
         \delta-\frac{p\delta T}{D}
    \end{pmatrix}.
    \end{equation}
Since $\lambda=\alpha+\beta \qi$ belongs to $\cD$, 
our assumption that $\delta_0(\lambda)>0$ implies the equality 
$\frac{T}{D}=\frac{T(\alpha,\beta,\delta_0(\lambda)}{D(\alpha,\beta,\delta_0(\lambda)}=\frac{1}{p}$.

We will proceed in 3 steps:

    In Step 1 we calculate and simplify derivatives of $\widehat{H}$.
   
    In Step 2 using row and column operation we reduce $3\times3$ determinant to $2\times2$ determinant.
    
    In Step 3 using Cauchy-Schwarz inequality we show that if $\mu$ is not a point mass, the $2\times 2$ determinant is non-zero.

\textbf{Step 1.} Before examining the required Jacobian determinant, let us simplify some of the entries of 
the Jacobian matrix of $\widehat{H}$.
We first calculate the derivative with respect to $\alpha$ of the first component of $\widehat{H}$:
		\begin{align*}
			&\frac{\partial \Re\left(\widehat{H}_{12}\left(\alpha,\beta,\delta\right)\right)}{\partial \alpha}
            =1-p\frac{\frac{\partial S}{\partial \alpha}}{D}+p\frac{S\frac{\partial D}{\partial \alpha}}{D^2}\\
            &\qquad\qquad\qquad=1-p\frac{T}{D}+2p\frac{\int_{\bR} \frac{(\alpha-x)^2}{\left((\alpha-x)^2+\beta^2+\delta^2\right)^2}d\mu(x)}{D}+p\frac{S\frac{\partial D}{\partial \alpha}}{D^2},
		\end{align*}
where in the last equality we used Proposition~\ref{Prop:STderivatives}. 
Since we observed that $1-p\frac{T}{D}=0$, we come to the conclusion that
		\begin{align*}
			\frac{\partial \Re\left(\widehat{H}_{12}\left(\alpha,\beta,\delta\right)\right)}{\partial \alpha}
            =2p\frac{\int_{\bR} \frac{(\alpha-x)^2}{\left((\alpha-x)^2+\beta^2+\delta^2\right)^2}d\mu(x)}{D}
            +p\frac{S\frac{\partial D}{\partial \alpha}}{D^2}.
		\end{align*}
We next calculate the derivative of the second component of $\widehat{H}$ with respect to $\beta$:	
\begin{align*}
			\frac{\partial \Im\left(\widehat{H}_{12}\left(\alpha,\beta,\delta\right)\right)}{\partial \beta}&=1+p\frac{(-T-\beta\frac{\partial T}{\partial \beta})D-\frac{\partial D}{\partial \beta}(1-\beta T)}{D^2}\\
            &=-p\frac{\beta \frac{\partial T}{\partial \beta}}{D}-p\frac{\frac{\partial D}{\partial \beta}(1-\beta T)}{D^2}.
		\end{align*}

Similarly, for the derivative of the third component of $\widehat{H}$ with respect to $\delta$ we get
		\begin{align*}
			\frac{\partial \Im \widehat{H}_{11}(\alpha,\beta,\delta)}{\partial \delta}=1-p\frac{T}{D}-p\delta\frac{\frac{\partial T}{\partial \delta} D-\frac{\partial D}{\partial \delta} T}{D^2}=-p\delta\frac{\frac{\partial T}{\partial \delta} D-\frac{\partial D}{\partial \delta} T}{D^2}.
		\end{align*}
		It is important to note that, based on Proposition~\ref{Prop:STderivatives}, an elementary 
        calculation gives that 
		\begin{align}\label{eq: Dderivatives}
			\frac{\partial D}{\partial \beta}-\frac{\beta}{\delta}\frac{\partial D}{\partial \delta}=-2T.
		\end{align}
        
		\textbf{Step 2.} We are now ready to process the required Jacobian determinant. Since the expression for every
        entry in that Jacobian matrix includes a factor of $p$, we can take out a $p^3$ from the determinant, 
        and we are left to show positivity for the determinant of the following matrix:
		\begin{align*}
			\begin{bmatrix}
				\frac{2}{D} \int_{\bR} \frac{(\alpha-x)^2}{\left((\alpha-x)^2+\beta^2+\delta^2\right)^2}d\mu(x)+\frac{S\frac{\partial D}{\partial \alpha}}{D^2}& -\frac{D\frac{\partial S}{\partial \beta}-S\frac{\partial D}{\partial \beta}}{D^2}&-\frac{D\frac{\partial S}{\partial \delta}-S\frac{\partial D}{\partial \delta}}{D^2}\\ \, \\
				\frac{-\beta \frac{\partial T}{\partial \alpha}D-(1-\beta T)\frac{\partial D}{\partial \alpha}}{D^2}& \frac{-\beta\frac{\partial T}{\partial \beta}D-\frac{\partial D}{\partial \beta}(1-\beta T)}{D^2}& \frac{-\beta\frac{\partial T}{\partial \delta}D-\frac{\partial D}{\partial \delta}(1-\beta T)}{D^2}\\\, \\
				-\delta \frac{\frac{\partial T}{\partial \alpha}D-T\frac{\partial D}{\partial \alpha}}{D^2} & -\delta \frac{\frac{\partial T}{\partial \beta}D-T\frac{\partial D}{\partial \beta}}{D^2}& -\delta\frac{\frac{\partial T}{\partial \delta} D-\frac{\partial D}{\partial \delta} T}{D^2}
			\end{bmatrix}
		\end{align*}
	Observe that the row operation $r_2:=r_2-\frac{\beta}{\delta} r_3$ gives
		\begin{align*}
			\begin{bmatrix}
				\frac{2}{D} \int_{\bR} \frac{(\alpha-x)^2}{\left((\alpha-x)^2+\beta^2+\delta^2\right)^2}d\mu(x)+\frac{S\frac{\partial D}{\partial \alpha}}{D^2}& -\frac{D\frac{\partial S}{\partial \beta}-S\frac{\partial D}{\partial \beta}}{D^2}&-\frac{D\frac{\partial S}{\partial \delta}-S\frac{\partial D}{\partial \delta}}{D^2}\\ \, \\
				-\frac{1}{D^2}\frac{\partial D}{\partial \alpha}& -\frac{1}{D^2}\frac{\partial D}{\partial \beta}& -\frac{1}{D^2}\frac{\partial D}{\partial \delta}\\\, \\
				-\delta \frac{\frac{\partial T}{\partial \alpha}D-T\frac{\partial D}{\partial \alpha}}{D^2} & -\delta \frac{\frac{\partial T}{\partial \beta}D-T\frac{\partial D}{\partial \beta}}{D^2}& -\delta\frac{\frac{\partial T}{\partial \delta} D-\frac{\partial D}{\partial \delta} T}{D^2}
			\end{bmatrix}.
		\end{align*}
		Next using second row we eliminate all derivatives of $D$ in the first and the third row, namely we perfom $r_1:=r_1+S r_2$ and $r_3:=r_3+\delta T r_2$, hence we obtain 
		\begin{align*}
			\begin{bmatrix}
				\frac{2}{D} \int_{\bR} \frac{(\alpha-x)^2}{\left((\alpha-x)^2+\beta^2+\delta^2\right)^2}d\mu(x)& -\frac{1}{D}\frac{\partial S}{\partial \beta}&-\frac{1}{D}\frac{\partial S}{\partial \delta}\\ \, \\
				-\frac{1}{D^2}\frac{\partial D}{\partial \alpha}& -\frac{1}{D^2}\frac{\partial D}{\partial \beta}& -\frac{1}{D^2}\frac{\partial D}{\partial \delta}\\\, \\
				-\frac{\delta}{D}\frac{\partial T}{\partial \alpha} & -\frac{\delta}{D}\frac{\partial T}{\partial \beta}& -\frac{\delta}{D}\frac{\partial T}{\partial \delta}
			\end{bmatrix}
		\end{align*}
		Next we will use relations between $\frac{\partial}{\partial \beta}$ and $\frac{\partial}{\partial \delta}$ applied to $S$ and $T$ found in Proposition~\ref{Prop:STderivatives} and similar properties when these two partial derivatives are applied to $D$ from \eqref{eq: Dderivatives}. Performing the column operation $c_2:=c_2-\frac{\beta}{\delta}c_3$ gives
		\begin{align*}
			\begin{bmatrix}
				\frac{2}{D} \int_{\bR} \frac{(\alpha-t)^2}{\left((\alpha-t)^2+\beta^2+\delta^2\right)^2}d\mu(t)& 0&-\frac{1}{D}\frac{\partial S}{\partial \delta}\\ \, \\
				-\frac{1}{D^2}\frac{\partial D}{\partial \alpha}& \frac{2T}{D^2}& -\frac{1}{D^2}\frac{\partial D}{\partial \delta}\\\, \\
				-\frac{\delta}{D}\frac{\partial T}{\partial \alpha} & 0& -\frac{\delta}{D}\frac{\partial T}{\partial \delta}
			\end{bmatrix}
		\end{align*}
    We expand determinant of the matrix above along the second column, clearly the (2,2) entry is positive. So in order to show that determinant is non-zero it suffices to show that
    \[
    -2 \frac{\partial T}{\partial \delta} \int_{\bR} \frac{(\alpha-t)^2}{\left((\alpha-t)^2+\beta^2+\delta^2\right)^2}d\mu(t)-\frac{\partial T}{\partial \alpha} \frac{\partial S}{\partial \delta}>0.
    \]
    
    \textbf{Step 3.}
    Denote $d\nu(t)=\frac{d\mu_x(t)}{((\alpha-t)^2+\beta^2+\delta^2)^2}$, then the inequality above becomes 
    \[
    -2 \frac{\partial T}{\partial \delta} \int_{\bR} (\alpha-t)^2d\nu(t)-\frac{\partial T}{\partial \alpha} \frac{\partial S}{\partial \delta}>0.
    \]
    Moreover from Proposition \ref{Prop:STderivatives} we have 
    \begin{align*}
        \frac{\partial T}{\partial\delta}&=-2 \delta \int_{\bR}d\nu(t),\\
        \frac{\partial T}{\partial \alpha}&=-2\int_{\bR}(\alpha-t)d\nu(t),\\
        \frac{\partial S}{\partial \delta}&=-2\delta \int_{\bR}(\alpha-t)d\nu(t),
    \end{align*}
    and the required inequality for determinant follows directly from Cauchy-Schwarz inequality in $L^2(\nu)$, by taking into account that $\nu$ cannot be a point mass.
\end{proof}

\subsection{The Jacobian outside the closure of the domain \boldmath{$\cD$}. }
\phantom{x}

\noindent
As we discussed under Assumption~\ref{assumption:16} mapping $H$ is well defined at matrices of the form $\begin{bmatrix}
    0&\lambda\\\overline\lambda&0
\end{bmatrix}$. In this subsection we show that the Jacobian $h$ at such point is positive, which shows that $H$ is locally invertible.

\begin{proposition}\label{prop:JacobianH-outside}
    With notations as in Proposition~\ref{prop:JacobianH}, 
    under Assumption \ref{assumption:16}, if $\lambda=\alpha+\qi\beta\notin \mathrm{cl}(\cD)$, 
    then we have
    \begin{equation}\label{eqn:Jacobian-positive-outside}
     \det\left(J_{\widehat{H}}(\alpha,\beta,0)\right)>0.
    \end{equation}
\end{proposition}
\begin{proof}
We recall that $\delta_0(\alpha+\qi \beta)=0$ outside the domain $\cD$. 

The Jacobian matrix of $\widehat{H}$ valued at $(\alpha,\beta,0)$ is of the form
\[ 
\begin{bmatrix}
   * & * & * \\
   * & *  & *\\
   0 &  0 & \frac{\partial \Im\left(\widehat{H}_{11}\left(\alpha,\beta,\delta\right)\right)}{\partial \delta}\Bigr|_{(\alpha,\beta,0)}
 \end{bmatrix}.
\]
The upper left $2\times 2$ matrix is the Jacobian matrix of the function $h$. In  Proposition~\ref{prop:outside-domain-injective} we showed that $h$ is analytic and injective on $\bC \setminus  (\mathrm{cl}(\cD)\cup\{\alpha_1, \alpha_2, \cdots, \alpha_k\})$ hence its Jacobian is nonzero. Moreover, 
Proposition~\ref{prop:12} and \eqref{eqn:limit-H11-over-delta} imply that 
\[
\frac{\partial \Im\left(\widehat{H}_{11}\left(\alpha,\beta,\delta\right)\right)}{\partial \delta}\Bigr|_{(\alpha,\beta,0)}=  \lim_{\delta\rightarrow{0}}\frac{\Im H_{11}(B(\lambda,\qi\delta))}{\delta}=1-\frac{p\cdot \varphi[(|\lambda-x|^2)^{-1}]}{|i+G_{\mu_x}(\lambda)|^2}\geq 0. 
\]
By Lemma \ref{lemma:subharmonic-1/F}, the function  $\lambda\mapsto\frac{\varphi[(|\lambda-x|^2)^{-1}]}{|i+G_{\mu_x}(\lambda)|^2}$ is subharmonic in $\mathbb{C}\setminus \spec(x)$. Hence, 
\[
\frac{\partial \Im\left(\widehat{H}_{11}\left(\alpha,\beta,\delta\right)\right)}{\partial \delta}\Bigr|_{(\alpha,\beta,0)}= 1-\frac{p\cdot \varphi[(|\lambda-x|^2)^{-1}]}{|i+G_{\mu_x}(\lambda)|^2}> 0
\]
for any $\lambda\notin (\mathrm{cl}(\cD)\cup \text{spec}(x))$.

If $\alpha\in\spec(x)$ is an atom of $\mu_x$ and $\alpha\notin \mathrm{cl}(\cD)$, then $\mu_x(\{\alpha\})>p$ by Proposition \ref{prop:atoms}. If follows from \eqref{eqn:D-over-T-general} that
\[
 \lim_{\delta\rightarrow{0}}\frac{\Im H_{11}(B(\lambda,\qi\delta))}{\delta}
  =\lim_{\delta\rightarrow{0}} \frac{T(\alpha,0,\delta)}{D(\alpha,0,\delta)}=\frac{1}{\mu_x(\{\alpha\})}.
\]
We have again 
\[
 \frac{\partial \Im\left(\widehat{H}_{11}\left(\alpha,\beta,\delta\right)\right)}{\partial \delta}\Bigr|_{(\alpha,\beta,0)}> 0
\]
for such $\lambda=\alpha$. 
We conclude that $\det\left(J_{\widehat{H}}(\alpha,\beta,0)\right)>0$.
\end{proof}

\subsection{Consequences for $h$}
\phantom{x}

\noindent
In this subsection we see some consequences of the positivity of the Jacobian of $H$.
\begin{corollary}
   \label{cor:Jh-positive-inside-domain}
    The Jacobian $J_h$ of $h$ is strictly positive inside the domain $\cD$.
\end{corollary}

\begin{proof}
    Let us denote $\widehat{H}(\alpha,\beta,\delta)=(\widehat{H}_1(\alpha,\beta,\delta),\widehat{H}_2(\alpha,\beta,\delta),\widehat{H}_3(\alpha,\beta,\delta))$. In this section we will view $h$ as a mapping from $\bR^2\to\bR^2$, and we have by the definition of $h$ that 
    \[h(\alpha,\beta)=(\widehat{H}_1(\alpha,\beta,\delta_0(\alpha,\beta)),\widehat{H}_2(\alpha,\beta,\delta_0(\alpha,\beta))).\]
    Direct calculation of the Jacobian matrix of $h$ gives
    \begin{align}\label{eqn:Jacobian_h}
        J_h=\begin{bmatrix}
            \frac{\partial H_1}{\partial \alpha}+\frac{\partial H_1}{\partial \delta} \frac{\partial \delta_0}{\partial \alpha}&\frac{\partial H_1}{\partial \beta}+\frac{\partial H_1}{\partial \delta} \frac{\partial \delta_0}{\partial \beta}\\
            \frac{\partial H_2}{\partial \alpha}+\frac{\partial H_2}{\partial \delta} \frac{\partial \delta_0}{\partial \alpha}&\frac{\partial H_2}{\partial \beta}+\frac{\partial H_2}{\partial \delta} \frac{\partial \delta_0}{\partial \beta}
        \end{bmatrix}
    \end{align}
    From the definition of $\delta_0$ we have
    \[H_3(\alpha,\beta,\delta_0(\alpha,\beta))=0.\]
    Differentiating this equation with respect to  $\alpha$ (for the first equation) and $\beta$ (for the second equation) gives
    \begin{align*}
        \frac{\partial H_3}{\partial \alpha}+\frac{\partial H_3}{\partial \delta}\frac{\partial \delta_0}{\partial \alpha}&=0,\\
        \frac{\partial H_3}{\partial \beta}+\frac{\partial H_3}{\partial \delta}\frac{\partial \delta_0}{\partial \beta}&=0.
    \end{align*}
    Hence we get
    \begin{align*}
        \frac{\partial \delta_0}{\partial \alpha}&=-\frac{\frac{\partial H_3}{\partial \alpha}}{\frac{\partial H_3}{\partial \delta}},\\
        \frac{\partial \delta_0}{\partial \beta}&=-\frac{\frac{\partial H_3}{\partial \beta}}{\frac{\partial H_3}{\partial \delta}}.
    \end{align*}
    Plugging this into Jacobian \eqref{eqn:Jacobian_h} and factoring out the denominator $\frac{\partial H_3}{\partial \delta}$ we get
    \begin{align*}
        J_h=\frac{1}{\frac{\partial H_3}{\partial \delta}}\begin{bmatrix}
            \frac{\partial H_1}{\partial \alpha}\frac{\partial H_3}{\partial \delta}-\frac{\partial H_1}{\partial \delta} \frac{\partial H_3}{\partial \alpha}&\frac{\partial H_1}{\partial \beta}\frac{\partial H_3}{\partial \delta}-\frac{\partial H_1}{\partial \delta} \frac{\partial H_3}{\partial \beta}\\
            \frac{\partial H_2}{\partial \alpha}\frac{\partial H_3}{\partial \delta}-\frac{\partial H_2}{\partial \delta} \frac{\partial H_3}{\partial \alpha}&\frac{\partial H_2}{\partial \beta}\frac{\partial H_3}{\partial \delta}-\frac{\partial H_2}{\partial \delta} \frac{\partial H_3}{\partial \beta}
        \end{bmatrix}
    \end{align*}
    A direct calculation of the determinant (note that the factor $\frac{1}{\frac{\partial H_3}{\partial \delta}}$ when calculating the determinant gets squared) of the above matrix shows that out of 8 summands two cancel, namely $\pm \frac{\partial H_1}{\partial \delta} \frac{\partial H_3}{\partial \alpha}\frac{\partial H_2}{\partial \delta} \frac{\partial H_3}{\partial \beta}$. The remaining 6 summands are directly related to the determinant of $3\times 3$ Jacobian, more precisely we have
    \begin{align*}
        \det(J_h)=\frac{1}{\frac{\partial H_3}{\partial \delta}} \det\begin{bmatrix}
            \frac{\partial H_1}{\partial \alpha}&\frac{\partial H_1}{\partial \beta}&\frac{\partial H_1}{\partial \delta}\\
            \frac{\partial H_2}{\partial \alpha}&\frac{\partial H_2}{\partial \beta}&\frac{\partial H_2}{\partial \delta}\\
            \frac{\partial H_3}{\partial \alpha}&\frac{\partial H_3}{\partial \beta}&\frac{\partial H_3}{\partial \delta}
        \end{bmatrix}.
    \end{align*}
    Since $\delta\mapsto\frac{H_3(\alpha,\beta,\delta)}{\delta}$ is monotone and non-negative for $\delta\geq \delta_0(\alpha,\beta)$, then  $\delta\mapsto H_3(\alpha,\beta,\delta)$ is monotone as a product of two monotone, non-negative functions. So $\frac{\partial H_3}{\partial \delta}$ is positive and the determinant of $3\times 3$ matrix is positive from Proposition~\ref{prop:JacobianH}.
\end{proof}

\begin{corollary}
$h^{-1}$ is differentiable on $\cM$.
\end{corollary}

\begin{proof}
This is an immediate consequence of Corollary \ref{cor:Jh-positive-inside-domain} and of the 
Inverse Function Theorem.
\end{proof}

Another consequence of the $3\times 3$ Jacobian being positive is injectivity of $h$, we first note that $h$ is continuous.

\begin{lemma}   \label{homeomorphism-lemma-A}
The function $h$ is continuous on $\mathbb{C}$.
\end{lemma}
\begin{proof} 
The continuity of $h$ follows from the continuity 
    of $\delta_0$ proved in Proposition \ref{prop:delta-continuity} and the fact that (under our running
    assumption $\spec(x)\subset\cD$) the functions used in the integrals defining $T$ and $S$ are bounded. 
\end{proof}

\begin{lemma}   \label{homeomorphism-lemma-B}
The
function $h$ is injective on $\bC$.
\end{lemma}

\begin{proof}
We already know from Proposition \ref{prop:outside-domain-injective} that $h$ is injective on $\bC\setminus\cD$.
    Thus it is enough to show that for any $\lambda_1\in\cD$ we have $h(\lambda_1)\neq h(\lambda)$ for any $\lambda\in\bC$ (so we look for an overlap of value of $h$ inside $\cD$ with any other point). 
	Towards contradiction assume that there is $\lambda_1=\alpha_1+\qi \beta_1\in\cD$ and $\lambda_2=\alpha_2+\qi\beta_2\in \bC$ such that $h(\lambda_1)=h(\lambda_2)$, then by the definition of $H$ we have
	\[H\left(
	\begin{bmatrix}
		\qi\delta_0(\lambda_1)&\lambda_1\\\overline{\lambda_1}&\qi\delta_0(\lambda_1)
	\end{bmatrix}\right)=\begin{bmatrix}
	0&h(\lambda_1)\\\overline{h(\lambda_1)}&0
	\end{bmatrix}
	\]
	and
	\[H\left(
	\begin{bmatrix}
		\qi\delta_0(\lambda_2)&\lambda_2\\\overline{\lambda_2}&\qi\delta_0(\lambda_2)
	\end{bmatrix}\right)=\begin{bmatrix}
		0&h(\lambda_2)\\\overline{h(\lambda_2)}&0
	\end{bmatrix}.
	\]
	So we have
	\[H\left(
	\begin{bmatrix}
		\qi\delta_0(\lambda_1)&\lambda_1\\\overline{\lambda_1}&\qi\delta_0(\lambda_1)
	\end{bmatrix}\right)=H\left(
	\begin{bmatrix}
	\qi\delta_0(\lambda_2)&\lambda_2\\\overline{\lambda_2}&\qi\delta_0(\lambda_2)
	\end{bmatrix}\right)\]
	In terms of $\widehat{H}$ from Proposition 5.2 we have
	$\widehat{H}(\alpha_1,\beta_1,\delta_0(\alpha_1+i\beta_1))=\widehat{H}(\alpha_2,\beta_2,\delta_0(\alpha_2+i\beta_2))$.
	
	Now we take two sequences $\delta_{1,n}$ and $\delta_{2,n}$ such that $\delta_{1,n}>\delta_0(\lambda_1)$ and $\delta_{1,n}\to\delta_0(\lambda_1)$ and similarly $\delta_{2,n}>\delta_0(\lambda_2)$ and $\delta_{2,n}\to\delta_0(\lambda_2)$. Which means that in the 3-dimensional space these two sequences converge to $(\alpha_1,\beta_1,\delta_0(\alpha_1+i\beta_1))$ and $(\alpha_2,\beta_2,\delta_0(\alpha_2+i\beta_2))$. From Proposition 5.2 we know that $H$ has local continuous inverse at $\widehat{H}(\alpha_1,\beta_1,\delta_0(\alpha_1+i\beta_1))$, because Jacobian determinant is non-zero at this point.
	
	We know that $\widehat{H}$ is injective on the set $\{(\alpha,\beta,\delta):\delta>\delta_0(\alpha+i\beta)\}$, as this is image of $\{(x,y,\varepsilon);\varepsilon>0\}$ under $\widehat\Omega$. Hence images under $\widehat H$ of $\{(\alpha_1,\beta_1,\delta_{1,n})\}$ and $\{(\alpha_2,\beta_2,\delta_{2,n})\}$ are disjoint and both define sequences converging to the same point $\widehat{H}(\alpha_1,\beta_1,\delta_0(\alpha_1+\qi \beta_1))=\widehat{H}(\alpha_2,\beta_2,\delta_0(\alpha_2+\qi \beta_2))$, this is because $\widehat{H}$ is continuous at any point of the form $(\alpha,\beta,\delta_0(\alpha+i\beta))$. Existence of two such sequences contradicts existence of local inverse of $\widehat{H}$ at $\widehat{H}(\alpha_1,\beta_1,\delta_0(\alpha_1+\qi \beta_1))$. More precisely we construct a sequence 
	\begin{align*}
		(a_n,b_n,\varepsilon_n)=\begin{cases}
			&\widehat{H}(\alpha_1,\beta_1,\delta_{1,n/2}) \mbox{ for $n$ even} 
			\\&\widehat{H}(\alpha_2,\beta_2,\delta_{2,(n+1)/2}) \mbox{ for $n$ odd}.
		\end{cases}
	\end{align*}
	For this sequence we have $(a_n,b_n,\varepsilon_n)\to \widehat{H}(\alpha_1,\beta_1,\delta_0(\alpha_1+\qi \beta_1))=\widehat{H}(\alpha_2,\beta_2,\delta_0(\alpha_2+\qi \beta_2))$.
	Since this sequence is contained in the domain of $\widehat{\Omega}$ we have 
	\begin{align*}
	\widehat{H}^{-1}(a_n,b_n,\varepsilon_n)=\Omega(a_n,b_n,\varepsilon_n)=\begin{cases}
	&(\alpha_1,\beta_1,\delta_{1,n/2}) \mbox{ for $n$ even} 
	\\&(\alpha_2,\beta_2,\delta_{2,(n+1)/2}) \mbox{ for $n$ odd}.
\end{cases}
\end{align*}
This contradicts existence of continuous inverse of $\widehat{H}$ at the point $\widehat{H}(\alpha_1,\beta_1,\delta_0(\alpha_1+\qi \beta_1))=\widehat{H}(\alpha_2,\beta_2,\delta_0(\alpha_2+\qi \beta_2))$.
	\end{proof}

\begin{corollary}\label{cor:62}
    Let $\cM=h(\cD)$ be the image of $\cD$ under the map $h$. The sets $\cD$ and $\cM$ are open subsets of $\mathbb{C}$. The function $h$ is a diffeomorphism from $\cD$ onto $\cM$.
\end{corollary}
\begin{proof}
 We infer from assertion $3^o$ in Proposition~\ref{prop:12} that $\cD$ is an open set in $\mathbb{C}$.
    From Lemma~\ref{homeomorphism-lemma-B} and Corollary~\ref{cor:h-differentiable-on-D} we know that $h$ is injective and differentiable on $\cD$. By Brouwer's invariance of domain theorem \cite[Corollary~19.8]{Bredon:1993}, it implies that $\cM$ is also open. Moreover, by Corollary~\ref{cor:Jh-positive-inside-domain}, the inverse of $h$ is also differentiable. Hence $h$ is a diffeomorphism from $\cD$ onto $\cM$.
\end{proof}

Next we show that, under Assumption \ref{assumption:16}, the function 
$h$ is a homeomorphism of the complex plane. 

\begin{theorem}\label{thm:h_is_homeomorphism}
    The function $h$ is a homeomorphism of $\bC$.
\end{theorem}
\begin{proof}
    It follows from general considerations that a function which is continuous, injective
and proper is a homeomorphism. 
We already know that $h$ is continuous and injective.
It was left to prove that for a compactly supported probability measure $\mu_x$ the function 
$h(\lambda)=\lambda-\frac{p}{\qi+G_{\mu_x}(\lambda)}$ is proper.

Pick a $K > 0$ such that $\mathrm{supp}(\mu_x)\subset [-K,K]$. For $\lambda \in \bC$ 
with $| \lambda | \geq 2K$ we then have $| \lambda - t | \geq | \lambda | /2$ for every 
$t \in \mathrm{supp} ( \mu_x )$, which implies that
\[
| G_{\mu_x}(\lambda) | \leq \int_\bR\frac{1}{|\lambda-t|}d\mu_x(t)\leq \frac{2}{ | \lambda |}.
\]
So if we also assume that $| \lambda | \geq 4$ then we get
\[
| \qi + G_{\mu_x}(\lambda)| \geq 1-|G_{\mu_x}(\lambda)| \geq 1 - \frac{2}{ | \lambda | } 
\geq \frac{1}{2},
\]
and hence that
\[
| h(\lambda) | \geq |\lambda|-\frac{p}{|\qi+G_{\mu_x}(\lambda)|}\geq |\lambda|-2p.
\]
The latter inequality implies that, for every $r > 0$, the disc $\{ z \in \bC \mid \, |z| \leq r \}$
has a bounded pre-image under $h$.  This, in turn, immediately implies that $h$ is proper. 
\end{proof}

\section{Formula for the density of the Brown measure}
\label{section:6}

For $\lambda = \alpha+\qi\beta\in\cD$, we let $z=s+\qi t=h(\lambda)\in \cM$, and thus have:
\begin{equation} \label{eqn:formula-h}
   \begin{cases}
    s=\Re\left(\widehat{H}_{12}\left(\alpha,\beta,\delta_0(\lambda)\right)\right)
                                       =\alpha-\frac{pS}{D}=\alpha-\frac{S}{T},   \\
    t=\Im\left(\widehat{H}_{12}\left(\alpha,\beta,\delta_0(\lambda)\right)\right)
                                       =\beta+p\frac{1-\beta T}{D}=\frac{1}{T},
  \end{cases}
\end{equation}
where we used the identity $\frac{T}{D}=\frac{1}{p}$ for $\lambda\in\cD$. 

Recall from the review in Section 2.1 that, in connection to our element of interest 
$a = x + \qi y$, we define 
\[
L_a(z)=\varphi(\log(|z-a|)) = \log \Delta(z-a)
\]
and 
\[
 L_a(z,\varepsilon)= \frac{1}{2} \varphi(\log(|z-a|^2+\varepsilon^2)).
\]

Following Lemma~\ref{lemma:22-Omega} and \eqref{eqn:H(B)}, we have 
 \begin{equation}
  \label{eqn:sub:6.2}
       \Omega\left( \begin{bmatrix}
                 \qi\varepsilon & z\\
                 \overline{z} & \qi\varepsilon            \end{bmatrix} \right)
                  =\begin{bmatrix}
                      \qi\delta(z,\varepsilon) & \lambda(z,\varepsilon)\\
                      \overline{\lambda(z,\varepsilon)} & \qi\delta(z,\varepsilon)
                  \end{bmatrix},
 \end{equation}
where we denote $\lambda= \lambda (z,\varepsilon)$ and $\delta=\delta(z,\varepsilon)$ 
to emphasize that $\lambda$ and $\delta$ are functions of 
$(z,\varepsilon)\in \mathbb{C}\times (0,\infty)$.

\begin{lemma}   \label{lemma:limits-Omega}
 Using the notations from \eqref{eqn:sub:6.2}, for any $z\in\mathbb{C}$, we have 
 \begin{enumerate}
     \item $\delta(z,\varepsilon)>\varepsilon$ for any $\varepsilon>0$ and $\lim_{\varepsilon\rightarrow{\infty}}\frac{\delta(z,\varepsilon)}{\varepsilon}=1$; and 
     \item $\lim_{\varepsilon\rightarrow{\infty}}\lambda(z,\varepsilon)=z-p\cdot\qi $.
 \end{enumerate}
 Moreover, the following limit exits
 \begin{equation}
  \label{eqn:integration-exists}
  \lim_{t\rightarrow{\infty}}\left(\int_1^{t}\delta(z,u) \varphi ( (|\lambda(z,u)-x|^2+\delta(z,u)^2)^{-1} -\frac{2}{1+u} \right)du
 \end{equation}
 and the limit is a real analytic function of $z$.
\end{lemma}
\begin{proof}
Recall that $H(\Omega(b))=b$ for any $b\in \bH^+(M_2(\mathbb{C}))$, and for $B=B(\alpha+\qi\beta,\qi\delta)$ we have 
\begin{align*}
    H(B)=B-\frac{p}{D}\begin{bmatrix}
  \qi\delta T&S-\qi(1-\beta T)\\
  S+\qi(1-\beta T)&\qi\delta T.
    \end{bmatrix},
\end{align*}
where  $D=D(\alpha,\beta,\delta)=\delta^2T^2+S^2+(1-\beta T)^2$.
From the definitions of $S$ and $T$ in \eqref{eqn:Sdefinition} and  \eqref{eqn:Tdefinition}, we can check that 
\[
  \lim_{\delta\rightarrow{\infty}}\frac{1}{D}\begin{bmatrix}
  \qi\delta T&S-\qi(1-\beta T)\\
  S+\qi(1-\beta T)&\qi\delta T.
    \end{bmatrix}=\begin{bmatrix}
         0 & -\qi\\
         \qi & 0
    \end{bmatrix}.
\]
Then (1) and (2) follow from Corollary~\ref{cor:imageSubord}.

By (2), we may assume that $\Vert \lambda(z,u)-x\Vert\leq M$ for some $M>0$. We then have 
\begin{align*}
          \delta(z,u) \varphi ( (|\lambda(z,u)-x|^2+\delta(z,u)^2)^{-1} )<\frac{1}{\delta(z,u)}<\frac{1}{u},
\end{align*}
and 
\[
     \frac{\delta(z,u)}{M^2+\delta(z,u)^2} <\delta(z,u) \varphi ( (|\lambda(z,u)-x|^2+\delta(z,u)^2)^{-1} ).
\]
Hence, 
\[
       \delta(z,u) \varphi ( (|\lambda(z,u)-x|^2+\delta(z,u)^2)^{-1} )>\frac{u}{M^2+u^2}
\]
provided that $\delta(z,u)> u>M$. Note that $\frac{2}{1+u}-\frac{u}{M^2+u^2}$ is integrable on $[1,\infty)$. It follows that the limit \eqref{eqn:integration-exists} exists. Moreover, the integrand in \eqref{eqn:integration-exists} is a real analytic function of $z$. We conclude that the limit is also a real analytic function of $z$.
\end{proof}

\begin{lemma}\label{lemma:extension-analyticity}
 For $z=s+\qi t\in\cM$, The map
 \[
   (s,t,\varepsilon)\mapsto L_{a}( s+\qi t,\varepsilon)
 \]
 extends to a real analytic function in some neighborhood of $(s, t, 0)$. Moreover,
 \[
   L_{a}(z,0)=\frac{1}{2}\int_0^{\infty}\left(\frac{2}{1+u} -\delta(z,u) \varphi ( (|\lambda(z,u)-x|^2+\delta(z,u)^2)^{-1} \right)du.
 \]
\end{lemma}
\begin{proof}
 Following Corollary~\ref{cor:imageSubord}, we have 
 \[
           \Omega\left( \begin{bmatrix}
                 \qi\varepsilon & z\\
                 \overline{z} & \qi\varepsilon            \end{bmatrix} \right)
                  =\begin{bmatrix}
                      \qi\delta(z,\varepsilon) & \lambda(z,\varepsilon)\\
                      \overline{\lambda(z,\varepsilon)} & \qi\delta(z,\varepsilon)
                  \end{bmatrix}.
\]
The subordination relation gives $G_{X+Y}\left(\begin{bmatrix}
                 \qi\varepsilon & z\\
                 \overline{z} & \qi\varepsilon            \end{bmatrix} \right)=G_X\left(\begin{bmatrix}
                      \qi\delta(z,\varepsilon) & \lambda(z,\varepsilon)\\
                      \overline{\lambda(z,\varepsilon)} & \qi\delta(z,\varepsilon)
                  \end{bmatrix} \right)$.
   By comparing $(1,1)$-entry, upon recalling \eqref{eq:entriesMatrixCauchy1}, we have
   \begin{align}
\label{eqn:6.2} 
    2 \cdot \frac{\partial L_{a}( z,\varepsilon)}{\partial\varepsilon}
       &=\varepsilon\varphi((|z- a|^2+\varepsilon^2)^{-1})\nonumber\\
       &=\Im \left( G_{X+Y, 11} \left(\begin{bmatrix}
                 i\varepsilon & z\\
                 \overline{z} & i\varepsilon            \end{bmatrix} \right)\right)\nonumber\\
      &=\Im \left( G_{X, 11} \left(
      \begin{bmatrix}
                      i\delta(z,\varepsilon) & \lambda(z,\varepsilon)\\
                      \overline{\lambda(z,\varepsilon)} & i\delta(z,\varepsilon)
                  \end{bmatrix} 
      \right)\right)\nonumber\\
       &=\delta(z,\varepsilon) \varphi ( (|\lambda(z,\varepsilon)-x|^2+\delta(z,\varepsilon)^2)^{-1} ). 
   \end{align} 
By integrating with respect to $\varepsilon$ from \eqref{eqn:6.2} we obtain that for $t > \varepsilon$ we have
  \begin{align*}
2(L_{a}(z,t)-L_{a}(z,\varepsilon))
      & = \int_\varepsilon^{t}
            u\varphi((|z-a|^2+u^2)^{-1}) d u\nonumber\\
      & =  \int_\varepsilon^{t}\delta(z,u) \varphi ( (|\lambda(z,u)-x|^2+\delta(z,u)^2)^{-1} )du.
  \end{align*}
On the left-hand side of above formula we are interested in estimating $L_a ( z, \varepsilon )$.
Towards that end, we will make $t \to \infty$ and use the fact that 
$\lim_{t\rightarrow\infty} \bigl( L_{a}(z,t)-\log (1+t) \bigr) =0$. 
Thus we rewrite the left-hand side of the above formula as
\begin{align*}
2(L_{a}(z,t)-L_{a}(z,\varepsilon))
& = -2 \,  L_{a}(z,\varepsilon)) + 2 \, \bigl( L_a (z,t) - \log (1+t) \bigr) + 2 \log (1+t)  \\
&= -2 \,  L_{a}(z,\varepsilon)) + 2 \, \bigl( L_a (z,t) - \log (1+t) \bigr) + 2 \int_0^t \frac{1}{1+u} \, du. 
\end{align*}

Putting this together with \eqref{eqn:integration-exists} 
in Lemma~\ref{lemma:limits-Omega} then gives us that 
\begin{align*}
    -2 L_{a}(z,0)&=\lim_{\varepsilon\rightarrow{0}}\int_\varepsilon^{1}\left(\delta(z,u) 
    \varphi ( (|\lambda(z,u)-x|^2+\delta(z,u)^2)^{-1} )-\frac{2}{1+u}\right)du\\
   &+ \lim_{t\rightarrow{\infty}}
      \int_1^{t}\left(\delta(z,u) 
      \varphi ( (|\lambda(z,u)-x|^2+\delta(z,u)^2)^{-1} )-\frac{2}{1+u}\right)du .
\end{align*}
Recall that $\lim_{u\rightarrow 0^+}\delta(z,u)=\delta_0 (h^{-1}(z))$ 
and $\lim_{u\rightarrow{0^+}}\lambda(z,u)=h^{-1}(z)$.
For $z\in h(\cD)$, by Proposition~\ref{prop:JacobianH}, the Jacobian matrix of the 
map $(s,t, \varepsilon)\mapsto (\alpha,\beta,\delta)$ has positive determinant in some neighborhood of $(s,t,0)$.  
 Therefore for $z\in h(\cD)$, by Proposition~\ref{prop:JacobianH} and Corollary~\ref{cor:imageSubord}, we deduce that the limit
\[
  \lim_{\varepsilon\rightarrow 0^+} \int_\varepsilon^{1}\delta(z,u) \varphi ( (|\lambda(z,u)-x|^2+\delta(z,u)^2)^{-1} )du
\]
exists and the limit is a real analytic function of $z$. This together with \eqref{eqn:integration-exists} guarantee that the map  $(s,t,\varepsilon)\mapsto L_{a}(s+i t,\varepsilon)$
 extends to a real analytic function in some neighborhood of $(x, y, 0)$. 
\end{proof}

\begin{proposition}
\label{prop:partial-derivative-0}
For $s_o + \qi t_o \in \cM$, the function $(s,t)\mapsto L_{a}(s+\qi t)$ is real analytic 
in some neighborhood of $(s_o, t_o)$. Moreover,
 \begin{equation}
   \frac{\partial L_{a}}{\partial{s}}
   (s_o + \qi t_o)  = \frac{\alpha_o - s_o}{t_o}
\ \mbox{ and }
\ \frac{\partial L_{a}}{\partial{t}} (s_o + \qi t_o) = \frac{\beta_o}{t_o},
 \end{equation}
where we denoted $h^{-1}(s_o+ \qi t_o) =: \alpha_o + \qi \beta_o$.
\end{proposition}

\begin{proof}
 The first assertion follows from Lemma~\ref{lemma:extension-analyticity}.
 By the subordination relation
   \[
     \Omega \left ( \begin{bmatrix}
       \qi\varepsilon & s+\qi t\\
       s-\qi t & \qi\varepsilon
\end{bmatrix}        \right) =\begin{bmatrix}
    \qi\delta_\varepsilon & \alpha_\varepsilon+\qi\beta_\varepsilon\\
    \alpha_\varepsilon-\qi\beta_\varepsilon & \qi\delta_\varepsilon
\end{bmatrix},
   \]
   and $\lim_{\varepsilon\rightarrow {0}}\delta_\varepsilon=\delta_0(x+\qi y)>0$, and 
    $\lim_{\varepsilon\rightarrow 0}(\alpha_\varepsilon+\qi\beta_\varepsilon)=\alpha+\qi\beta$. Denote $\lambda_\varepsilon=\alpha_\varepsilon+\qi \beta_\varepsilon$,
    recall that the subordination relation implies that
     \begin{align*}
      &\varphi((z-x-\qi y)^*(|z-x-\qi y|^2+\varepsilon^2)^{-1})
   =\varphi( (\lambda_\varepsilon-x)^*(|\lambda_\varepsilon-x|^2+\delta_\varepsilon^2)^{-1} ).
     \end{align*}
For $z=s+\qi t$, by the regularity result in Lemma~\ref{lemma:extension-analyticity}, 
we are allowed to flip the derivative and the limit in the following  below. Using \eqref{eq:entriesMatrixCauchy1}, we get 
 \begin{align*}
  \frac{\partial L_{a}(s+\qi t)}{\partial{z}}
  &=\frac{\partial}{\partial {z}}\left(\lim_{\varepsilon\rightarrow {0}}L_{a}(z,\varepsilon)\right)\\
  &=\lim_{\varepsilon\rightarrow {0}}\left(\frac{\partial}{\partial {z}} L_{a}(z,\varepsilon)\right)\\
   &=\lim_{\varepsilon\rightarrow {0}} \frac{1}{2} \varphi((z-a)^*(|z-a|^2+\varepsilon^2)^{-1})\\
   &=\lim_{\varepsilon\rightarrow {0}} 
      \frac{1}{2} \varphi( (\lambda_\varepsilon-x)^*(|\lambda_\varepsilon-x|^2+\delta_\varepsilon^2)^{-1} )\\
   &= \frac{1}{2} \varphi( (\lambda-x)^*(|\lambda-x|^2+\delta_0(\lambda)^2)^{-1} ),
 \end{align*}
 where $\lambda=\alpha+\qi\beta=h^{-1}(s+\qi t)$ and $\delta_0(\lambda)>0$, and similarly 
  \begin{align*}
  \frac{\partial L_{a}(s+\qi t)}{\partial{\overline{z}}}
    =\varphi( (\lambda-x)(|\lambda-x|^2+\delta_0(\lambda)^2)^{-1} ).
 \end{align*}
 Hence, using Equation~\eqref{eqn:formula-h} we obtain 
 \begin{align*}
  \frac{\partial L_{a}}{\partial{s}} (s_o + \qi t_o)
    =S=(\alpha_o - s_o)T=\frac{\alpha_o -s_o}{t_o},
 \end{align*}
 and
 \begin{align*}
  \frac{\partial L_{a}}{\partial{t}} (s_o + \qi t_o)
   =\beta T=\frac{\beta_o}{t_o}.
 \end{align*}
\end{proof}

\begin{theorem} \label{thm:63}
   The Brown measure $\mu_a$ of $a=x+\qi y$ is absolutely continuous on $\cM$. Its density $f(s,t)$ for any $s+\qi t\in \cM$ is given by 
\begin{equation}  
        f(s,t)=\frac{1}{2\pi}\left[\frac{1}{t}\left(\frac{\partial \alpha}{\partial s}
           +\frac{\partial \beta}{\partial t}\right)-\frac{1}{t}-\frac{\beta}{t^2}\right],
    \end{equation}
where $\alpha+\qi\beta=h^{-1}(s+\qi t)$.
\end{theorem}
\begin{proof} Since (as seen in Lemma \ref{lemma:extension-analyticity}) $L_a$ is real analytic in a neighborhood of $(s,t,0)$, 
it follows that $f(s,t)$ can be retrieved as 
$\frac{1}{2 \pi } \left(\frac{\partial^2}{\partial s^2}+ \frac{\partial^2}{\partial t^2}\right)L_{a}( s+\qi t)$.
But, in view of Proposition~\ref{prop:partial-derivative-0}, we have 
    \begin{align*}
\left(\frac{\partial^2}{\partial s^2}+ \frac{\partial^2}{\partial t^2}\right)L_{a}( s+\qi t)
         & = \frac{\partial}{\partial s}\frac{\alpha-s}{t} +\frac{\partial}{\partial t}\frac{\beta}{t}\\
         & = \frac{1}{t}\left(\frac{\partial \alpha}{\partial s}
             +\frac{\partial \beta}{\partial t}\right)-\frac{1}{t}-\frac{\beta}{t^2}.
    \end{align*}
\end{proof}

\begin{proposition}  \label{prop:65}
If Assumption~\ref{assumption:16} holds,  using notations in Proposition \ref{prop:outside-domain-injective},
we have $\mu_{a} (h(\{\mathbb{C}\backslash\left(\mathrm{cl}(\cD)\cup\{\alpha_1,\ldots,\alpha_k\}\right)\}))=0$.
\end{proposition}

\begin{proof}
Let $\lambda=\alpha+\qi\beta\in \mathbb{C}\,\backslash\left(\mathrm{cl}(\cD)\cup\{\alpha_1,\ldots,\alpha_k\}\right)$ 
and $z=h(\lambda)$.  Recall that, by Proposition \ref{prop:JacobianH-outside}, the determinant
     $\det\left(J_{\widehat{H}}(\alpha,\beta,0)\right)>0$ and the map $\widehat{H}$ is locally invertible in some neighborhood of $(\alpha,\beta, 0)$. Note that $\lambda\notin\text{spec}(x)$ and hence $\lambda-x$ is invertible. 
     We apply the same argument as for Lemma~\ref{lemma:extension-analyticity}, and deduce that the map
 \[
   (s,t,\varepsilon)\mapsto L_a (s+\qi t,\varepsilon)
 \]
 extends to a real analytic function in some neighborhood of $(s, t, 0)$, where $s+\qi t=h(\alpha+\qi \beta)$. 
Consequently,  the map $(s,t)\mapsto L_a (s+\qi t)$ is a real analytic function
in some neighborhood of $(x,y)$. 

Hence we can exchange the order of taking derivative and the limit. We have 
\begin{align*}
\frac{\partial L_a}{\partial z} (s + \qi t)
& =\frac{\partial}{\partial {z}}\left(\lim_{\varepsilon\rightarrow {0}}L_a ( z,\varepsilon)\right)\\
& =\lim_{\varepsilon\rightarrow {0}}\left(\frac{\partial}{\partial {z}} L_a (z,\varepsilon)\right)\\
& =\lim_{\varepsilon\rightarrow {0}}\frac{1}{2}\varphi((z-a)^*(|z-a|^2+\varepsilon^2)^{-1})\\
& =\lim_{\varepsilon\rightarrow {0}}
     \frac{1}{2}\varphi( (\lambda_\varepsilon-x)^*(|\lambda_\varepsilon-x|^2+\delta_\varepsilon^2)^{-1} )\\
&=\frac{1}{2}\varphi\left( (\lambda-x)^{-1} \right),
\end{align*}
 where $\lambda=\alpha+\qi\beta=h^{-1}(s+\qi t)$ and 
  \begin{align*}
  \frac{\partial L_a}{\partial{\overline{z}}} (s + \qi t )
    =\frac{1}{2}\varphi\left( (\overline{\lambda}-x)^{-1} \right).
 \end{align*}
Moreover, by Proposition~\ref{prop:outside-domain-injective}, the function $h$ is holomorphic 
at $\lambda$ and locally invertible at $z=h(\lambda)=\lambda+\frac{p}{-\qi-G_{\mu_x}(\lambda)}$. 
Hence, $h^{-1}$ is holomorphic at $z$ and 
 \[
   \frac{\partial}{\partial\overline{z}}\varphi((h^{-1}(z)-x)^{-1})=0,
\ \mbox{ that is, } 
\ \frac{\partial^2 L_a}{\partial{\overline{z}}\partial{z}} (s+ \qi t) =0.
 \]
Therefore $z =s + \qi t$ does not belong to $\text{supp}(\mu_a)$.
\end{proof}

\vspace{10pt}

\section{The atoms of the Brown measure $\mu_a$}

Throughout this section, we suppose that $0<p<1$ in which case $\mu_y(\{0\})=1-p>0$. 
If $t$ is an atom of $\mu_x$ such that $\mu_x(\{t\})<p$, it is shown in the proof 
for (1) of Proposition \ref{prop:43} that $t\in \cD$. Suppose now $\mu_x(\{t\})>p$;
then, motivated by \cite[Theorem 7.4]{BV1998regularity} for self-adjoint model, one expects 
the Brown measure of $a = x+\qi y$ to have an atom at $t$ with weight $\mu_x(\{t\})-p$. 
The main result of this section is the following. 

\begin{theorem}
 \label{thm:71}
Suppose that $0 < p < 1$, and let $t_o \in \mathbb{C}\setminus\mathrm{cl}(\cD)$ be an isolated atom of 
$x$ with $\mu_x(\{t_o\})>p$. Under Assumption \ref{assumption:16}, we then have 
$h(t_o)=t_o\in\mathbb{C}\setminus\mathrm{cl}(\cM)$ and 
  \[
    \mu_{a}(\{t_o \})= \mu_x(\{t_o\})-p.
  \]
\end{theorem}

\begin{lemma}
 \label{lemma:72}
 Let $z_o\in\mathbb{C}$ and suppose that $L_a(z)$ 
 is harmonic in $\{ z \in \bC : 0 < |z - z_o| < r \}$.
 Then, $z_o$ is an isolated atom of $\mu_a$ and moreover $\mu_a(\{z_o\})$ can be calculated by
 \[
   \mu_a(\{z_0\})=\frac{1}{\pi\qi} \oint_{|z-z_o|=r_o}\frac{\partial L_a(z)}{\partial z}dz,
 \]
 for any $0<r_o<r$.
\end{lemma}

\begin{proof}
Fix an $r_0$ such that $0 < r_o <r$, and denote $D = \{z:|z-z_o|\leq r_0\}$.
Recall that in distributional sense we have 
\[
\mu_a= \frac{2}{\pi} \frac{\partial}{\partial \overline{z}}\frac{\partial}{\partial {z}} L_a(z).
\]
By applying Green's formula for $f(z)=\frac{\partial}{\partial {z}} L_a(z)$, we have 
    \[
\oint_{\partial D} \frac{\partial}{\partial {z}} L_a(z)\,dz
=
2i \iint_D \frac{\partial}{\partial \overline{z}}\frac{\partial}{\partial {z}} L_a(z)\,dA,
\]
where $dA$ denotes the standard area measure.

Since $L_a$ is harmonic in $D \setminus\{z_0\}$, we have
\[
\frac{\partial}{\partial \overline{z}}
\frac{\partial}{\partial z}L_a(z)=0
\]
in the classical sense on $D\setminus\{z_o\}$. Hence, the distribution
$\frac{\partial}{\partial \overline{z}}
\frac{\partial}{\partial z}L_a(z)$
is supported at $z_o$.
The result then follows by substituting the integral back into the distributional definition of $\mu_a$.
\end{proof}

\begin{proof}[Proof of Theorem \ref{thm:71}]
We first note that $G_{\mu_x}(t_o)=\infty$ and hence by Proposition \ref{prop:outside-domain-injective} we have $h(t_o)=t_o$. 
We recall that $h(\mathbb{C}\setminus\mathrm{cl}(\cD))=\mathbb{C}\setminus\mathrm{cl}(\cM)$.
By repeating a calculation done in the proof of Proposition \ref{prop:65}, we find that for any 
$\lambda\in \mathbb{C}\setminus(\mathrm{cl}(\cD)\cup\spec(x))$ we have 
\begin{align}\label{eqn:701}
\frac{\partial L_{a}}{\partial{z}} (s + \qi t )
  = \frac{1}{2} \varphi((\lambda-x)^{-1})= \frac{1}{2} G_{\mu_x}(\lambda),
\end{align}
where $z = s + \qi t$ is related to $\lambda$ via the formula 
\begin{equation}
  \label{eqn:71}
     z=h(\lambda)=\lambda-\frac{p}{\qi+G_{\mu_x}(\lambda)}.
\end{equation}
When we solve for $G_{\mu_x}(\lambda)$ in \eqref{eqn:71} and substitute the resulting expression into 
\eqref{eqn:701}, it follows that for $\lambda\in \mathbb{C}\setminus(\mathrm{cl}(\cD)\cup\spec(x))$
\begin{equation}   \label{eqn:73}
  \frac{\partial L_{x+\qi y}(s+\qi t)}{\partial{z}}= \frac{1}{2} \frac{p}{\lambda-z}- \frac{\qi}{2}
  =\frac{1}{2}\frac{p}{h^{-1}(z)-z}- \frac{\qi}{2},
\end{equation}
since $h$ is a homeomorphism. 
We have $h^{-1}(t)=t$, thus $\frac{\partial L_{x+\qi y}(s+\qi t)}{\partial{z}}$ has a singularity at $z=t$. 
We need to say that $h^{-1}$ is holomorphic, and hence $L_{x+\qi y}$ is harmonic $\frac{\partial }{\partial{\overline{z}}}\frac{\partial }{\partial{z}}L_{x+\qi y}=0$.

By applying Lemma \ref{lemma:72} and then by using Equation (\ref{eqn:73}), we find that 
\begin{align*}
    \mu_a(\{t_o\})&=\frac{1}{\pi\qi} \oint_{|z-t_o|=r_o}\frac{\partial L_a(z)}{\partial z}dz\\
    &=\frac{1}{2\pi\qi}\oint_{|z-t_o|=r_0}\left(\frac{p}{h^{-1}(z)-z}-\qi\right)dz,
\end{align*}
which will allow us to compute $\mu_a ( \{ t_o \} )$ as residue of an analytic function.
One can check directly that $h'(t_o)=1-\frac{p}{\mu_x(\{t_o\})}$. Hence we have 
\[
(h^{-1})'(t_o) = \frac{1}{ h' ( \, h^{-1} (t_o) \, )} = \frac{1}{h'(t_o)} = \frac{\mu_x(\{t\})}{\mu_x(\{t_o\})-p}.
\]
It then follows that
\[
 \lim_{z\rightarrow t}(z-t)\times \left(\frac{p}{h^{-1}(z)-z}-\qi\right)=\mu_x(\{t_o\})-p.
\]
Therefore, this shows that $\mu_{a}(\{t_o\})=\mu_x(\{t_o\})-p$.
\end{proof}

\vspace{0.5cm}

\begin{proof}
[Proof of Theorem \ref{thm:18}]

$\ $

\noindent

The claim follows by combining Theorem \ref{thm:63}, Theorem \ref{thm:71}, and Proposition \ref{prop:65}.

Theorem \ref{thm:63} shows that the Brown measure $\mu_a$ is absolutely continuous on $\cM$ and gives its density there.

Theorem \ref{thm:71} identifies the possible atoms of $\mu_a$ as $\alpha_1,\alpha_2,\ldots,\alpha_k$, and shows that
\[
\mu_a(\{\alpha_i\})=\mu_x(\{\alpha_i\})-p,
\qquad i=1,2,\ldots,k.
\]

Finally, Proposition \ref{prop:65} implies that
\[
\mu_a\!\left(h\!\left(\mathbb{C}\setminus
\bigl(\mathrm{cl}(\cD)\cup\{\alpha_1,\alpha_2,\ldots,\alpha_k\}\bigr)\right)\right)=0.
\]
Since $h$ is a homeomorphism, $h(\cD)=\cM$, and $h(\alpha_i)=\alpha_i$ for
$i=1,\ldots,k$, we have
\[
h\!\left(\mathbb{C}\setminus
\bigl(\mathrm{cl}(\cD)\cup\{\alpha_1,\alpha_2,\ldots,\alpha_k\}\bigr)\right)
=
\mathbb{C}\setminus
\bigl(\mathrm{cl}(\cM)\cup\{\alpha_1,\alpha_2,\ldots,\alpha_k\}\bigr).
\]
Therefore,
\[
\mu_a\!\left(\mathbb{C}\setminus
\bigl(\mathrm{cl}(\cM)\cup\{\alpha_1,\alpha_2,\ldots,\alpha_k\}\bigr)\right)=0,
\]
which proves that the support of $\mu_a$ is contained in
$\mathrm{cl}(\cM)\cup\{\alpha_1,\alpha_2,\ldots,\alpha_k\}$.
\end{proof}

\noindent 

{\bf AI Declaration} ChatGPT 5.6 Sol was used to examine and verify preliminary arguments concerning the discrete part of the Brown measure. These arguments contributed to the development of the proof of Theorem 7.1 and of the part of Proposition 3.9 concerning the injectivity of $h$ at atoms. The final proofs were written and independently verified by the authors, who take full responsibility for their correctness.

\noindent

{\bf Acknowledgement.} We are grateful to Serban Belinschi and to Brian Hall for 
insightful comments concerning the first version of this paper.
We would like to thank the anonymous referee for the very careful reading of the manuscript and for many valuable comments and suggestions. We are especially grateful for pointing out some issues in an earlier version of the paper and for suggestions that led to the present formulation of Theorem \ref{thm:18}, which now also describes the atomic part of the Brown measure.

\bibliographystyle{amsplain}

\bibliography{brown}

\end{document}